\documentclass[10pt]{article}
\usepackage[top=25truemm,bottom=25truemm,left=30truemm,right=30truemm]{geometry}
\usepackage{authblk}

\usepackage{graphicx} % Required for inserting images
\usepackage{color} % Required for inserting images
\usepackage[
	colorlinks=true, 
	linkcolor=black, 
	citecolor=blue, 
	urlcolor=blue]{hyperref}
\usepackage{amsmath,amssymb,amsthm,mathtools}
\usepackage[capitalize,noabbrev]{cleveref}
\usepackage{mysty,mythm}

\usepackage{lineno}
\makeatletter
	
	\@addtoreset{equation}{section}
\makeatother

\theoremstyle{definition}

\usepackage{framed} % frame
\usepackage{txfonts} % math font
\usepackage{xfrac}
\usepackage{enumitem}
\usepackage{subfigure}
\usepackage{bm}
\usepackage{natbib}
\bibpunct[:]{(}{)}{,}{a}{}{,} % (author, 2022) のようにする
\usepackage{algorithm}
\usepackage{algpseudocode}
\usepackage{empheq}
\usepackage{multicol}
\usepackage{overpic}
\usepackage{extarrows}
\usepackage{here}
\usepackage{comment}
\allowdisplaybreaks

\usepackage{physics}
\title{A market-based efficient matching mechanism for crowdsourced delivery systems with demand/supply elasticities}
\author[1]{Yuki Oyama\thanks{Corresponding Author: oyama@bin.t.u-tokyo.ac.jp}} 
\affil[1]{Department of Civil Engineering, The University of Tokyo, Tokyo, Japan}

\author[2]{Takashi Akamatsu\thanks{Corresponding Author: akamatsu@plan.civil.tohoku.ac.jp}}
\affil[2]{Graduate School of Information Sciences, Tohoku University, Miyagi, Japan}

\date{\today}

\begin{document}
\maketitle

\begin{abstract}
Crowdsourced delivery (CSD), or crowd-shipping, is an emerging business model that leverages the underutilized or excess capacity of individual drivers to fulfill delivery tasks. 
% By utilizing existing trips of ordinary drivers, CSD is expected to not only perform fast and cost-effective deliveries but also reduce the negative environmental impact associated with the use of dedicated delivery vehicles.
This paper presents a general formulation of a larege-scale two-sided CSD matching problem, considering demand/supply elasticity, heterogeneous preferences of both shippers and drivers, and task-bundling. We propose a set of methodologies to solve this problem. First, we reveal that the fluid-particle decomposition approach of \cite{akamatsu_oyama_csd_2024} can be extended to our general formulation. This approach decomposes the original large-scale matching problem into a fluidly-approximated task partition problem (master problem) and small-scale particle matching problems (sub-problems). We propose to introduce a truthful and efficient auction mechanism to the sub-problems, which enables the observation of privately perceived costs for each shipper/driver. Furthermore, by finding a theoretical link between the auction problems and random/parturbed utility theory, we succeed in accurately reflecting the information collected from the auctions to the master problem in a theoretically consistent manner. This reduces the master problem to a smooth convex optimization problem, theoretically guaranteeing the computational efficiency and solution accuracy of the fluid approximation. 
% The task partition problem also achieves global matching efficiency by evaluating as its objective the optimal value functions of the sub-problems based on the perturbed utility theory. 
Second, we transform the master problem into a traffic assignment problem (TAP) based on a task-chain network representation. This transformation overcomes the difficulty in enumerating task bundles. Finally, we formulate the dual problem of the TAP, whose decision variable is only a price/reward pattern at market equilibrium, and develop an efficient accelerated gradient descent method. The numerical experiments clarify that the proposed approach drastically reduces the computational cost of the general and complex CSD matching problem ($\sim$ 700 times faster than a naive method) without sacrificing accuracy of the optimal solution (mostly within 0.5\% errors).

\hfill\break%
\noindent\textit{Keywords}: Crowdsourced delivery; Matching market design; Fluid–particle decomposition; Auction mechanism; Perturbed utility theory; Traffic assignment problem
\end{abstract}

\newpage
\tableofcontents % 目次を生成
\newpage

\section{Introduction}
% The widespread diffusion of e-commerce and rising customer expectations of delivery 
The rapid growth of e-commerce and rising customer expectations of fast, low-price, and punctual delivery have significantly increased the demand for last-mile parcel delivery. Last-mile delivery is challenging for traditional couriers with a limited number of drivers due to the large volume of small but frequent delivery requests, necessitating novel and cost-effective delivery solutions. 
In recent years, the widespread adoption of smart mobile and communication technology has stimulated the development of on-demand mobility services and the gig economy. Crowdsourced delivery (CSD), or crowd-shipping, is an emerging business model for city logistics inspired by this trend. CSD leverages the underutilized or excess capacity of individual drivers, such as commuters, travelers, and shoppers, to fulfill delivery tasks \citep{le2019supply}. 
%## CSD basic concept, challenges
% Crowdsourced delivery (CSD) fulfills delivery tasks by utilizing the excess capacity of ordinary drivers' private vehicles that may not be fully used. As a promising solution to urban logistics, 
By utilizing existing trips of ordinary drivers, CSD is expected to not only perform fast and cost-effective deliveries but also reduce the negative environmental impact associated with the use of dedicated delivery vehicles, such as traffic, pollutant emissions, and energy consumption \citep{mladenow2016crowd, paloheimo2016transport, simoni2020potential}.

\begin{figure}[t]
    \centering
    \includegraphics[width=0.90\textwidth]{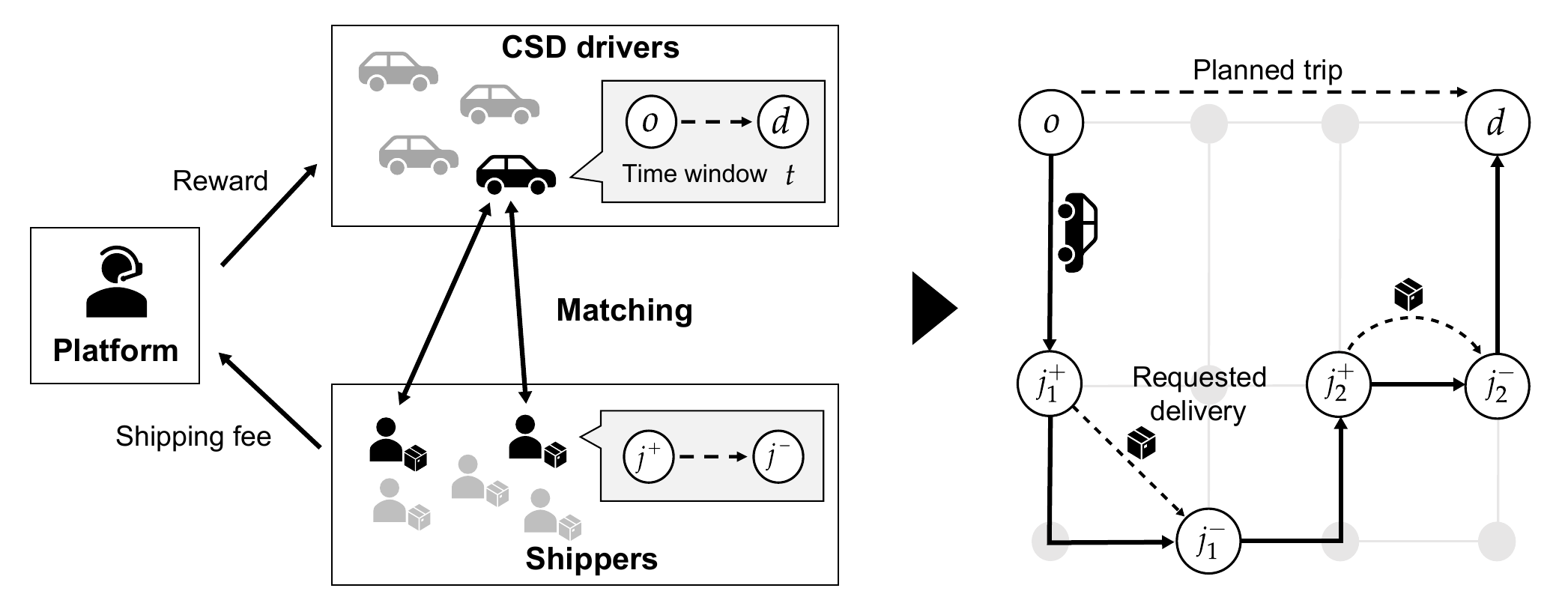}
    \caption{Conceptual diagram of the CSD system considered in this study.}
    \label{fig:CSD}
\end{figure}

%matching "market design"
In this paper, we study a two-sided CSD matching market of \cref{fig:CSD}, where the platform acts as an intermediary to match shippers (demand) and drivers (supply). Drivers are ordinary individuals who have their own origin-destination pairs and perform delivery tasks on the way by making detours.
% This study deals with a CSD system in which ordinary drivers who travel for their own purposes perform delivery tasks on the way from the origin to the destination by making detours.
% the platform matches delivery tasks received from individual shippers with ordinary drivers who travel for their own purposes and perform delivery tasks on the way from origin to destination by making detours. 
% This is often called \textit{en-route matching} and has been extensively studied in the literature \citep{alnaggar2021crowdsourced}. Specifically, we study a two-sided CSD matching market of \cref{fig:CSD}, where the platform plays as an intermediary to match shippers (demand) and drivers (supply). 
Both shippers and drivers decide on their participation in the CSD system based on the prices/rewards offered by the platform, and these decisions vary among individuals.
% We consider stochastic participation and choice behavior of both shippers and drivers as functions of the prices/rewards offered by the platform. 
We also consider task-bundling, that is the possibility of drivers performing multiple delivery tasks during a single trip. Although realistic, these settings are associated with the issues of (i) computational complexity, (ii) supply and demand uncertainties, and (iii) heterogeneity in individual preferences, and make the matching problem very challenging in today's e-commerce context with a large demand for parcel delivery. We aim to simultaneously solve all three issues, while they have been only independently addressed in the literature, as explained below. 
% that there is no existing study that simultaneously addresses all three issues, and this study fills the gap.

%# complexity
First, the CSD matching problem is a large-scale combinatorial optimization problem that is generally difficult to solve in a reasonable computational time. This problem becomes even harder when task-bundling is considered, as the number of possible task bundles explodes as that of tasks increases.
To reduce computational complexity, some studies decompose the problem into smaller-scale and tractable matching problems, e.g., by consolidating pickup points to a limited number of depots \citep{wang2016towards, archetti2016vehicle, arslan2019crowdsourced} or clustering tasks based on their geographical locations \citep{huang2021solving, elsokkary2023crowdsourced, simoni2023crowdsourced}. However, a two-stage optimization that relies on arbitrary problem decomposition does not ensure global matching efficiency. It may also generate a partitioned set of delivery tasks with a wide range of distances, particularly in cases where the pickup and delivery locations of requests are widely distributed across the entire network, leading to an inefficient matching pattern. In addition, most existing studies assume one-to-one matching in which only a single task is assigned to each driver \citep[e.g.,][]{archetti2016vehicle, soto2017matching, dayarian2020crowdshipping, ccinar2024role, akamatsu_oyama_csd_2024}, while assigning a bundle of tasks to a driver can significantly reduce the total transportation cost and the compensation to be paid due to the sub-additivity of the costs \citep{gansterer2018centralized, mancini2022bundle, wang2023joint}. 
Considering task-bundling adds extra complexity to the matching problem and has never been addressed in a large-scale setting in the literature.
% This approach can also hedge against the supply uncertainty \citep{yang2024freelance}. However, allowing drivers to perform multiple tasks adds further complexity to the matching problem, making it even harder to solve.

Second, CSD systems often consider ordinary individuals, such as commuters and travelers, as supply resources. These drivers are called occasional/ad-hoc/freelance drivers \citep{archetti2016vehicle, arslan2019crowdsourced, yang2024freelance}, and they are not dedicated drivers but travel anyway for their purposes. Therefore, the participation of drivers in the CSD system is stochastic and depends on compensation \citep{le2019modeling, dayarian2020crowdshipping}, i.e., the supply is price-elastic. % and their preferences may be heterogeneous
The decision of drivers not to work for CSD leads to supply uncertainty and can significantly influence matching efficiency and operational costs \citep{hou2022optimization, yang2024freelance}. To deal with supply uncertainty in CSD matching, recent studies model drivers' acceptance behavior of dispatching orders as a function of compensation \citep{dahle2019pickup, le2021designing, barbosa2023data, ccinar2024role, hou2022optimization, hou2025reinforced}. Integrating a driver acceptance behavior model, \cite{hou2022optimization, hou2025reinforced}, and \cite{ccinar2024role} design personalized rewards through optimization and achieve a high acceptance rate and cost savings. However, these studies test their approaches only on small-scale problems and assume a one-to-one matching in which only a single task is offered to each driver. 
Moreover, most studies ignore the impact of reward values on demand, i.e., demand elasticity, despite the fact that increased rewards (salary) may lead to increased shipping fees and thus demand uncertainty. %except for \cite{le2021designing},

Finally, the willingness of drivers to work for the CSD and the willingness of shippers to pay for delivery can be very heterogeneous between individuals \citep{le2019influencing, le2019supply, dayarian2020crowdshipping}. Finding an optimal matching pattern requires information on these heterogeneous preferences, but they are generally unobservable \citep{akamatsu_oyama_csd_2024}. This unobserved heterogeneity generates stochasticity of behavior and supply/demand uncertainty. Collecting choice data is a way to address heterogeneity \citep{le2019modeling, barbosa2023data, hou2022optimization}, but such data do not necessarily elicit true preferences because the truly preferred alternatives for individuals may not be available/offered in the observed choice situations.
In contrast, several previous works study market-based CSD matching with bidding \citep{kafle2017design, allahviranloo2019dynamic, su2023exact, mancini2022bundle}. Performing auctions with a \textit{truthful} market mechanism, such as the Vickrey-Clarke-Groves (VCG) mechanism \citep{vickrey1961counterspeculation, clarke1971multipart, groves1973incentives}, enables the observation of individual preferences through bids. However, the VCG mechanism is computationally expensive and does not apply to a large-scale problem, and previous studies on bid-based CSD systems do not present mechanisms to elicit true preferences. 
% The Vickrey-Clarke-Groves (VCG) mechanism \citep{vickrey1961counterspeculation, clarke1971multipart, groves1973incentives} is one of the market mechanisms that ensure truth telling and efficiency, and performing auctions with the VCG mechanism is a way to observe individual preferences. However, the VCG mechanism does not apply to a large-scale problem as it needs to solve an optimization problem repeatedly. 

In our previous study \citep{akamatsu_oyama_csd_2024}, we proposed a fluid-decomposition (FPD) approach to deal with a large-scale CSD matching problem with unobserved heterogeneous driver preferences. The approach decomposes the matching problem into task partition (master problem) and individual task-driver matching within smaller groups of drivers (sub-problems). We applied the VCG mechanism to each sub-problem to observe true perceived costs through bids. Moreover, the FPD approach performs task partition so that the global matching efficiency is achieved. To approximate the optimal value functions of the sub-problems without explicitly solving them, we exploit the random utility modeling (RUM) framework and fluidly approximate the master problem. This way, the FPD approach addresses the issues of (i) computational complexity and (ii) heterogeneous preferences and achieves a highly accurate matching pattern with a very efficient computation. However, the previous study solves one-to-one matching and ignores demand and supply elasticities.

% In summary, solving the CSD matching problem involves three main issues: (i) demand/supply elasticities, (ii) preference heterogeneity and unobservability, and (iii) computational complexity. To the best of our knowledge, no existing study addresses the above three issues simultaneously. This study fills this gap. Specifically, we extend the work of \cite{akamatsu_oyama_csd_2024} who present the market-based efficient matching mechanism for CSD systems by proposing a fluid-particle decomposition (FPD) approach. The FPD approach decomposes the matching problem into task partition (master problem) and individual task-driver matching within smaller groups of drivers (sub-problems). By exploiting the random utility framework, the fluidly approximated master problem efficiently partitions delivery tasks while achieving global matching efficiency with high accuracy. In addition, the VCG auction mechanism is applied to the small-scale sub-problems so that the true perceived utilities of drivers are elicited. 
%As such, the mechanism simultaneously addresses the challenges of (i) and (ii) mentioned above.  

In this paper, we study a more general and complex matching problem that takes into account demand and supply elasticities and task-bundling, which are practical requirements for CSD systems. We model a two-sided market: On one side, shippers decide the time window to request delivery tasks, and on the other side, drivers decide which delivery tasks to perform.
We also consider opt-out options so that shippers and drivers can leave the system when prices/rewards do not meet their preferences. In addition, drivers are allowed to perform multiple delivery tasks during a single trip.

\begin{figure}[t]
    \centering
    \includegraphics[width=\textwidth]{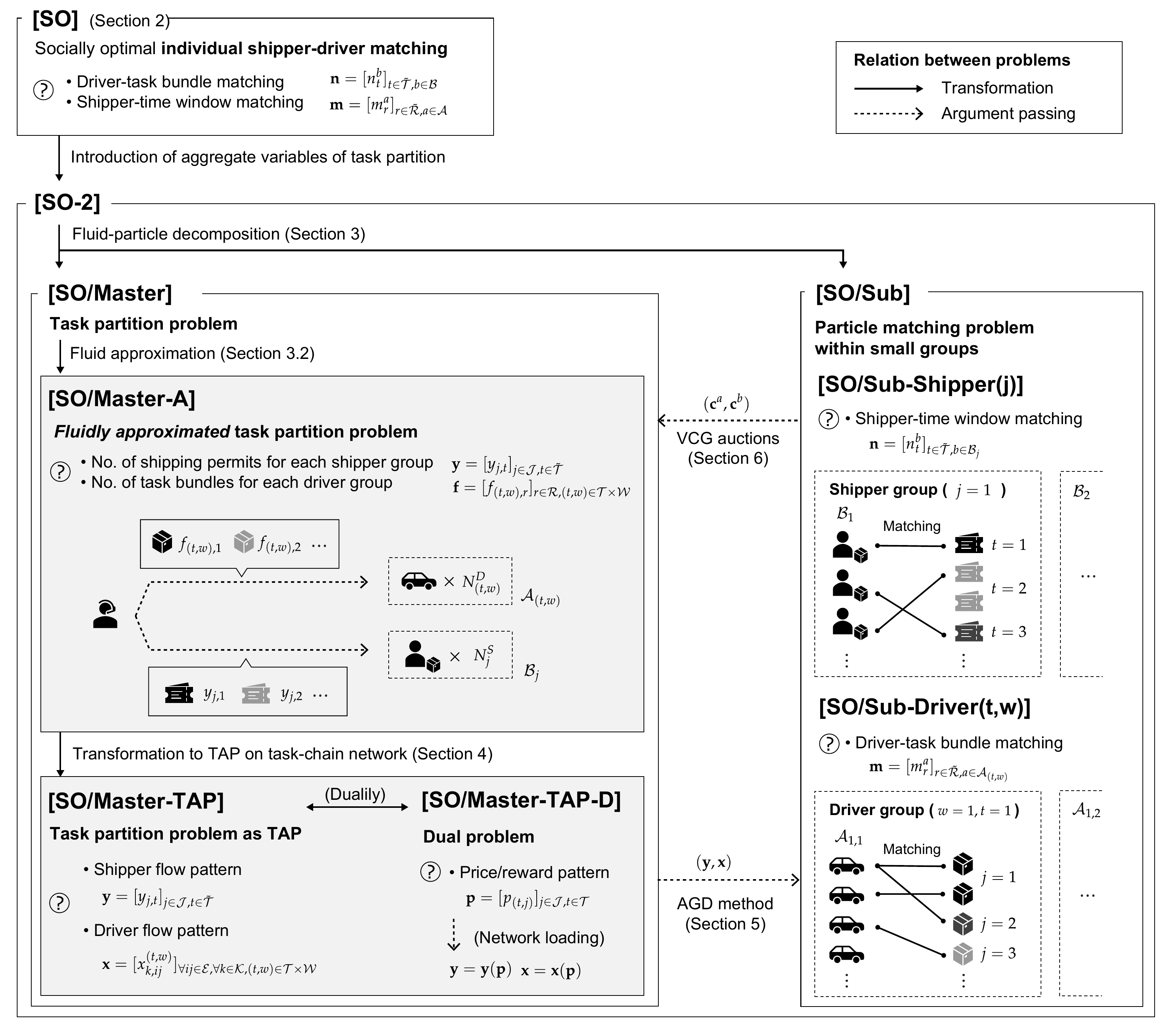}
    \caption{The full picture of the approach of this study.}
    \label{fig:flow}
\end{figure}

To solve this more generalized matching problem, we propose novel and efficient algorithms by extending the FPD approach. \cref{fig:flow} shows a complete picture of the approach of this study. The FPD approach decomposes the original matching problem into a fluidly approximated task partition problem (master problem) and particle matching problems (sub-problems) within driver and shipper subgroups. 
%Considering demand and supply elasticities, our master (task partition) problem simultaneously 
Instead of finding an individual matching pattern $(\mathbf{m}, \mathbf{n})$ all at once, the master problem decides the number of shipping permits for each shipper group ($\mathbf{y}$) and the number of task bundles for each driver group ($\mathbf{f}$). Furthermore, to overcome the prohibitive complexity of enumerating task bundles, we transform the master problem into a traffic assignment problem (TAP) on a task-chain network, which reduces to finding an edge flow pattern $\mathbf{x}$ instead of $\mathbf{f}$. We further exploit its dual formulation, whose decision variable is only a price/reward pattern $\mathbf{p}$, and develop a fast gradient-based solution algorithm. These decomposition and transformation techniques dramatically reduce the number of decision variables and computational costs. The sub-problems are small enough for us to quickly solve and apply the VCG mechanism that enables the observation of individual preferences through bids. As such, our proposed approach simultaneously addresses the issues of (i) computational complexity, (ii) supply and demand elasticity, and (iii) individual heterogeneity, which is, to the best of our knowledge, the first achievement in the literature.
% Prices and rewards are optimized so that market equilibrium is achieved.
% Moreover, we relax the assumptions of \cite{akamatsu_oyama_csd_2024} on the number of tasks for each driver. The assignment of multiple tasks to a driver adds further complexity to the matching problem.
% To efficiently solve the problem, this study represents the task-driver matching based on a task-chain network and casts the task partition problem as a capacitated traffic assignment problem (CTAP). We then develop a dual-type algorithm to solve the CTAP, which significantly reduces the number of decision variables and efficiently performs task partition. 

The notable contributions of this study are summarized as follows:
% \begin{itemize}
%     \item In \cref{sec:problem}, we define and formulate a socially optimal matching problem for a two-sided CSD market. We explicitly consider demand/supply elasticity depending on the price/compensation by modeling the stochastic choice behavior of both shippers and drivers, including opt-out options. We also perform the theoretical
%     \item In \cref{sec:mechanism}, we 
% \end{itemize}

\renewcommand{\theenumi}{\arabic{enumi}}
\renewcommand{\labelenumi}{(\theenumi)}
\begin{enumerate}
    \item \textbf{Formulation}: We formulate a socially optimal CSD matching problem, in which we consider demand/supply elasticity, heterogeneous preferences of both shippers and drivers, and task-bundling, and that the number of shippers and drivers can be very large. This formulation is general and simultaneously considers the three issues (i)-(iii) explained above.
    We also theoretically prove that the optimality condition coincides with the market equilibrium condition. Matching, pricing, and compensation are simultaneously optimized under the market equilibrium. 
    % We also prove that the optimality condition coincides with the market equilibrium condition.
    % \item We explicitly consider demand/supply elasticity depending on the price/compensation. We model the choice behavior of both shippers and drivers, including opt-out options.
    % \item We deal with the assignment of a bundle of tasks to a driver. 
    \item \textbf{Approach}: We reveal that the FPD approach of \cite{akamatsu_oyama_csd_2024} can be extended for our generalized formulation. The extended approach hierarchically decomposes the matching problem into a master problem for task partition and sub-problems for particle matching. The sub-problems for shipper/driver groups are small enough to quickly solve and apply a truthful and efficient auction mechanism, which ensures the observation of true perceived costs. 
    The master problem is fluidly approximated and exploits the perturbed utility theory to evaluate the optimal value functions of the sub-problems in a theoretically consistent manner.
    \item \textbf{Algorithms}: For the general and complex CSD matching problem, we develop efficient algorithms without sacrificing the accuracy of the optimal solution.
    To overcome the difficulty in enumerating task bundles, we introduce a task-chain network representation on which we cast the task partition problem as a TAP problem. Moreover, by exploiting its dual formulation, we develop an efficient algorithm based on the accelerated gradient descent method. The proposed transformation and algorithm drastically reduce the number of decision variables compared to the original problem, and their efficiency is demonstrated through numerical experiments. 
    % We also use the network representation and present a scalable edge-based VCG mechanism to solve particle matching problems.
\end{enumerate}

% The remainder of the paper is structured as follows. \ldots
% \Cref{sec:review} provides the literature review focusing on CSD with demand/supply elasticity. 

%
\section{Matching problem of crowdsourced delivery}\label{sec:problem}
% This section defines and formulates the matching problem for the CSD system. We formulate the matching problem as a combinatorial optimization problem to achieve the socially optimal state, and show that its optimality conditions can be interpreted as market equilibrium conditions. The major challenges in solving the matching problem are also discussed.

This study models a two-sided CSD market in which a CSD platform acts as an intermediary between a set of shippers and a set of drivers aiming to achieve \textit{socially optimal matching}. While most previous studies regard the platform as a profit maximizer, studying socially optimal matching is essential for public entities to appropriately implement policies and regulations for CSD markets to avoid ``market failure'' \citep{akamatsu_oyama_csd_2024}. In this section, we first formulate the socially optimal matching problem and then show that its optimality conditions coincides with market equilibrium. Note that the list of notation that frequently appears in this study is provided in \Cref{app:notation}.

\subsection{Problem definition and assumptions}
Let $\cG\equiv (\cN,\cL)$ be a directed graph representing a transportation network where $\cN$ and $\cL$ are the sets of nodes and links, respectively, and each link $kl\in\cL$ is associated with the link travel time $t_{kl}$. This study assumes that the travel time is static and deterministic and that the travels of CSD participants do not affect it.
Also, let $\cT \equiv \{1,2,\ldots,T\}$ be the set of time windows. Each time window is assumed to be long enough (e.g., a few hours) for drivers to perform their own trips or delivery tasks. Upon this setting, we define and assume the main components of the CSD system as follows.

%\equiv\{1,2,\dots,N\}
% Then, we consider a CSD system. The main components of the system are defined and assumed as follows.

\begin{itemize}
    \item \textbf{Delivery task}. A delivery task is a unit of delivery operations and must be delivered from pickup and delivery locations on network $\cG$. We assume that tasks are distinguished by only their pickup-delivery location pairs, and denote a task category by $j = (j^+, j^-) \in \cJ \subseteq\cN\times\cN$. In other words, all delivery tasks for the same pickup-delivery location pair $(j^+, j^-)$ are homogeneous. We assume each task must be completed within a time window and never done across multiple windows. 
    %The set of tasks to be delivered during window $t$ is defined as $\cJ$.
    % We assume that all delivery tasks for the same pickup-delivery pair during the same window are homogeneous and that a task is never performed across multiple time windows. 
    %This can be realized by grouping delivery operations into delivery-task units so that every delivery task is as homogeneous as possible. 
    %The number of tasks delivered between $(r,s)$ is denoted by $n_{rs}$.
    \item \textbf{Platform}. The CSD platform is an intermediary and aims to efficiently match shippers (demand) and drivers (supply). We assume that the platform manager is a non-profitable agent, and the matching objective is the maximization of social surplus.
    To this end, the platform opens two-sided markets for both shippers and drivers. On one side, the platform issues \textit{shipping permits} for each time window for each group of shippers, and in the market shippers submit bids to obtain the preferred permit.
    % task $j$ and each window $t$ and opens a market to sell them to the group of shippers $\cB_{j}$ who want to deliver goods for $j$.
    On the other side, the platform assigns the delivery tasks received from shippers to drivers, based on the bids submitted by CSD drivers.   
    % matches the tasks received from shippers with drivers through the market for drivers. Each CSD driver submits bids based on their private cost for performing a bundle of tasks, based on which the assignment and reward for each task are determined. 
    
    % receives delivery requests from shippers and outsources the delivery tasks to CSD drivers. The platform aims to efficiently match shippers and drivers so as to maximize the social surplus, which is the sum of the shippers' and drivers' perceived utilities. The platform manager is assumed to be a non-profitable agent. In other words, all revenue from shippers is assumed to be transferred to drivers through compensation. 
    % The manager employs a dedicated driver to perform a task for $(r,s)$ if no CSD driver is matched with the task. We assume that the operation cost of a dedicated vehicle is $\bar{c}_{rs}$ and consider it to be the upper bound of the reward $w_{rs}$ for outsourcing. 
    
    % \item \textbf{Two-sided market}. 
    % The platform first partitions the tasks to be delivered within time window $t$ into each
    % The winner of each auction performs the delivery task and is rewarded $w_{rs}$.
    
    \item \textbf{Shipper}. A shipper is an agent who has a unit parcel that must be delivered between pickup and delivery points. We assume each shipper has only one delivery task of a single unit, and therefore, we categorize shippers by the task's pickup-delivery pair $j = (j^+, j^-)$. We define the set of shippers $\cB_{j}$ who may request tasks for $j$ and its population $N^S_j = |\cB_{j}|$. 
    Shippers are assumed to have differenct preferences for time windows in which their parcels are delivered. As such, a shipper decides whether and in which time window $t$ to ship their parcel. The choice set is $\tilde{\cT} = \cT \cup \{0\}$ where $t=0$ represents the opt-out option (i.e., canceling the use of the CSD service). %delivering the task themself
    % , and this set is common for all shippers. 
    % or deliver it themself, and if outsourcing it, then chooses a time window $t \in \cT$. 
    % This decision of a shipper $b \in \cB_{j}$ is made so as to minimize their perceived cost $c^b_t$ plus the outsourcing fee $p_{(t,j)}$. 
    % If a shipper $b \in \cB_{j}$ outsources the delivery task during time window $t$, their perceived cost is $c^b_t$ plus the outsourcing fee $p_{(t,j)}$ decided by the platform. Otherwise, the perceived cost is $c^b_0$, mainly representing the transportation cost for delivering the task themself. We assume $c^b_0 - c^b_t > 0, \, \forall t \in \cT, \, \forall b \in \cB_j$. 
    
    \item \textbf{CSD driver}. A CSD driver is an ordinary network traveler who plans to travel for their purpose between an origin-destination (OD) pair $w = (o,d) \in \cW \subseteq \cN\times\cN$ during a time window $t \in \cT$, and they can perform delivery tasks on the way to their destination. For instance, an origin and a destination can be a residential zone and a central business district. We define the set of drivers $\cA_{(t,w)}$ traveling between OD pair $w$ during time window $t$ and its population $N^D_{(t,w)} = |\cA_{(t,w)}|$.
    % Let $a \in \cA_{(t,w)}$ denote the set of drivers traveling during the time window $t$ between an origin-destination (OD) pair $w = (o,d) \in \cW \subseteq \cN\times\cN$. 
    We assume that each driver can perform $K$ tasks at maximum but only one task can be executed at a time (i.e., the next task can be executed after the one has been completed). Therefore, each driver's decision is whether and which \textit{task bundle} $r$ to perform considering the detour cost and rewards for the tasks. The choice set is defined as $\tilde{\cR} = \cR\cup\{\emptyset\}$ where $\cR$ is the set of task bundles and $r = \emptyset$ represents not performing any task.
    % Therefore, the utility of a driver depends on a bundle of tasks to perform. 
    %, while the disutility of the driver not performing a delivery task is the perceived cost $c^a_{(o,d)}$ to travel between $(o,d)$ 
\end{itemize}

% Note that the single unit delivery per driver is a strong assumption and ignores the possiblity of a driver performing tasks with multiple times of pick-ups \citep[as in][]{archetti2016vehicle, soto2017matching,
% dayarian2020crowdshipping}, but this is more applicable for en-route matching \citep{alnaggar2021crowdsourced}. 
% The simplicity also allows for a clear presentation of our theoretical framework, which is the main focus of this paper. However, the fluid-particle decomposition approach of this paper can be extended to the case with multiple tasks per driver and/or price elasticity, and we leave the extension for future study (also see \Cref{sec:conclusion}).

\subsection{Behavior of agents}
\subsubsection{Shipper's behavior}
The behavior of a shipper $b \in \cB_j$ is to decide whether and in which time window $t$ to ship a parcel for pickup-delivery location pair $j$, to minimize their perceived cost.
Let $n^b_t \in \{0,1\}$, $\forall t \in \tilde{\cT} \equiv \cT \cup \{0\}$ be the choice variable of shipper $b \in \cB_{j}$. Here $n^b_0 = 1$ indicates that shipper $b$ does not use the CSD service with a perceived cost $c^b_0$ (e.g., representing the cost of requesting a professional courier or deliverying the parcel themself), while $n^b_t = 1, \forall t \in \cT$ indicates that they ship the parcel in time window $t$ with a perceived cost $c^b_t$ plus the shipping fee $p_{(t,j)}$. Hence, the optimal behavior of the shipper $b \in \cB_j$ is represented as the following minimization problem:
\begin{align}
    % [\text{SO}]\qquad
    \min_{\mathbf{n}_{b(j)}}\quad&n^b_0 c^b_0 + \sum_{t\in\cT}n^b_t(c^b_t + p_{(t,j)}) =
    \sum_{t \in \tilde{\cT}} n^b_t c^b_t + \sum_{t\in\cT}n^b_t p_{(t,j)}\label{eq:ship_obj}\\
    \subto\quad&\sum_{t \in \tilde{\cT}} n^b_t = 1 \qquad \label{eq:indv_shipper_consv}\\
    &n^b_t \in \{0,1\} \qquad \forall t\in \tilde{\cT},\label{eq:n_binary}
\end{align}
where the constraint (\ref{eq:indv_shipper_consv}) represents a single choice of the shipper.

\subsubsection{Driver's behavior}
In this study, a driver can perform at maximum $K$ tasks during their trip, and we represent a task bundle as a list $r = [j_1, \ldots, j_k]$ with $j_{1:k} \in \cJ$ and $k \leq K$.
The behavior of a driver $a \in \cA_{(t,w)}$ who travels between OD pair $w$ in time window $t$ is to decide whether or not they performs any delivery task and then choose a task bundle $r$ from the bundle set $\cR$ to minimize their perceived detour cost $c^a_r$ minus the reward $p_{(t,r)}$ associated with bundle $r$, wherein $p_{(t,r)}$ is defined as the sum of the rewards of the tasks included in $r$. Therefore, the entire choice set is $\tilde{\cR} = \cR\cup\{\emptyset\}$ where $r = \emptyset$ represents not performing any task. 
% Note that we define $p_{(t,r)}$ as the sum of the rewards of the tasks included in $r$, that is,
% \begin{align}\label{eq:bundle_price}
%     p_{(t,r)} = \sum_{j \in \cJ}\delta_{r,j} p_{(t,j)},
% \end{align}
% where $\delta_{r,j} \in \mathbb{Z}_{+}$ is the number of tasks for $j$ included in $r$.

Let $m^a_r \in \{0,1\}, \, \forall r \in \tilde{\cR}$ be the choice variable of driver $a \in \cA_{(t,w)}$, and then their optimal behavior is represented as the following minimization problem:
\begin{align}
    % [\text{SO}]\qquad
    \min_{\mathbf{m}_{a(t,w)}}\quad&m^a_{\emptyset} c^a_{\emptyset} + \sum_{r\in \cR} m^a_r (c^a_r - p_{(t,r)}) = \sum_{r\in \tilde{\cR}} m^a_r c^a_r - \sum_{r\in \cR} m^a_r \sum_{j \in \cJ}\delta_{r,j} p_{(t,j)}
    %\sum_{j \in r} p_{(t,j)}
    \label{eq:drv_obj}\\
    \subto\quad&\sum_{r\in \tilde{\cR}} m^a_r = 1
    \label{eq:indv_drv_consv}\\
    &m^a_r \in \{0,1\} \qquad \forall r \in \tilde{\cR},\label{eq:m_binary}
\end{align}
where $\delta_{r,j} \in \mathbb{Z}_{+}$ is the number of tasks for $j$ included in $r$. %the reward $p_{(t,r)}$ associated with bundle $r$ is the sum of the rewards of the tasks included in $r$, and

\subsection{System optimal matching problem}
Next, we formulate the socially optimal (SO) matching problem based on the optimization problems mentioned above. To do so, we consider the following demand-supply constraint:
\begin{align}
    \label{eq:demand_supply}
    \sum_{b\in\cB_{j}}n^b_t \leq \sum_{w \in \cW}\sum_{a\in \cA_{(t,w)}}\sum_{r\in \cR} m^a_r \delta_{r,j}\qquad \forall j \in \cJ, \,  \forall t \in \cT
\end{align}
which ensures that the number of drivers who are willing to deliver tasks $j$ in time window $t$ is equal to or greater than the number of shippers who desire to ship their parcels between $j$ in $t$.

The social cost to be minimized is the sum of shippers' and drivers' perceived costs \eqref{eq:ship_obj} and \eqref{eq:drv_obj}. By the assumption of a nonprofit platform, the following equation representing a money transfer between shippers and drivers holds:
\begin{align}
    \label{eq:money_transf}
    \sum_{j\in\cJ}\sum_{b\in\cB_{j}}\qty(\sum_{t\in\cT}n^b_t p_{(t,j)}) =\sum_{t\in\cT}\sum_{w \in \cW}\sum_{a\in \cA_{(t,w)}}\qty(\sum_{r\in \cR} m^a_r \sum_{j \in \cJ}\delta_{r,j} p_{(t,j)}) 
    %&= \sum_{t\in\cT}\sum_{w \in \cW}\qty(\sum_{a\in \cA_{(t,w)}}\sum_{r\in \cR} m^a_r \sum_{j \in r} p_{(t,j)}) \nonumber\\
\end{align}
which cancels out the money-related terms in \eqref{eq:ship_obj} and \eqref{eq:drv_obj}. 
This means that we only have to consider the perceived costs of shippers and drivers (without the fees/rewards) in the SO matching problem:
% As a result, the optimization problem [SO] for achieving the system optimal matching of delivery tasks and CSD drivers is formulated as follows.
\begin{align}
    [\text{SO}]\qquad
    \min_{\mathbf{m},\mathbf{n}}\quad&
    z(\mathbf{m},\mathbf{n}) \equiv \sum_{j\in\cJ}\qty(\sum_{b\in\cB_{j}} \sum_{t\in\tilde{\cT}} n^b_t c^b_t) +
    \sum_{t\in\cT}\sum_{w \in \cW}\qty(\sum_{a\in\cA_{(t,w)}} \sum_{r\in \tilde{\cR}} m^a_r c^a_r)
    \label{eq:SO_obj}\\
    \subto\quad& (\ref{eq:demand_supply}),\nonumber\\
    &\sum_{t\in\tilde{\cT}} n^b_t = 1 \qquad \forall b \in \cB_{j}, \, \forall j \in \cJ, \label{eq:indv_shipper_consv_all}\\
    &n^b_t\in\{0,1\} \qquad \forall t \in\tilde{\cT}, \, \forall b \in \cB_{j}, \, \forall j \in \cJ, \label{eq:n_binary_all}\\
    &\sum_{r\in \tilde{\cR}} m^a_r = 1 \qquad \forall a \in \cA_{(t,w)}, \, \forall (t,w) \in \cT \times \cW
    \label{eq:indiv_driver_consv_all}\\
    &m^a_r \in \{0,1\} \qquad \forall r \in \tilde{\cR}, \, \forall a \in \cA_{(t,w)}, \, \forall (t,w) \in \cT \times \cW\label{eq:m_binary_all}
\end{align}
% where the first constraint \eqref{eq:indv_driver_consv} is the condition that each driver performs one delivery task, and the second constraint \eqref{eq:indiv_task_consv} states that the supply of drivers must be equal to or less than the number of delivery tasks. % at most

\subsection{Optimality conditions and market equilibrium}
Due to the 0--1 binary constraints \eqref{eq:n_binary_all} and \eqref{eq:m_binary_all}, the problem [SO] is a combinatorial optimization problem. % that is hard to solve in general
However, because the constraint matrices are totally unimodular, we can relax the problem [SO] as linear programming by replacing the binary constraints \eqref{eq:n_binary_all} and \eqref{eq:m_binary_all} with non-negative constraints $n^b_t\geq 0$ and $m^a_r\geq 0$. Based on the linearly relaxed formulation, we can analyze the optimal state of the matching problem [SO] as follows. 

We define the Lagrangian $\cL$ as
\begin{align}
    \cL(\mathbf{m},\mathbf{n},\bm{\mu},\bm{\rho},\bm{\lambda}) 
    &\coloneqq z(\mathbf{m},\mathbf{n}) 
    +\sum_{j\in\cJ}\sum_{b\in\cB_{j}}\mu_{b(j)} \qty(1 - \sum_{t \in \tilde{\cT}}n^b_t)  + \sum_{t\in\cT}\sum_{w \in \cW}\sum_{a\in\cA_{(t,w)}} 
    \rho_{a(t,w)} \qty(1 - \sum_{r\in \tilde{\cR}}m^a_r) \nonumber\\   
    &~~~~~ +\sum_{t\in\cT}\sum_{j\in \cJ}\lambda_{(t,j)}\qty(\sum_{b\in\cB_{j}}n^b_t -\sum_{w \in \cW}\sum_{a\in \cA_{(t,w)}}\sum_{r\in \cR} m^a_r \delta_{r,j})
\end{align}
where $\bm{\mu} \geq \bm{0}$, $\bm{\rho} \geq \bm{0}$ and $\bm{\lambda} \geq \bm{0}$ are the Lagrangian multipliers associated with the constraints \eqref{eq:indv_shipper_consv_all}, \eqref{eq:indiv_driver_consv_all} and \eqref{eq:demand_supply}, respectively. The optimality conditions of [SO] are
\begin{align}\label{eq:opt_cond_sop}
    \begin{dcases}
        \pdv{\cL}{n^b_t} = c^b_t + \lambda_{(t,j)} - \mu_{b(j)} \geq 0 ,\, n^b_t \geq 0 ,\, n^b_t \pdv{\cL}{n^b_t} = 0, \\
        \pdv{\cL}{m^a_r} = c^a_r - \sum_{j\in\cJ}\delta_{r,j}\lambda_{(t,j)} + \rho_{a(t,w)} \geq 0 ,\, m^a_r \geq 0 ,\, m^a_r \pdv{\cL}{m^a_r} = 0, \\
        \pdv{\cL}{\lambda_{(t,j)}} = \sum_{b\in\cB_{j}}n^b_t -\sum_{w \in \cW}\sum_{a\in \cA_{(t,w)}}\sum_{r\in \cR} m^a_r \delta_{r,j} \geq 0,\, \lambda_{(t,j)} \geq 0,\, \lambda_{(t,j)}\qty(\sum_{b\in\cB_{j}}n^b_t -\sum_{w \in \cW}\sum_{a\in \cA_{(t,w)}}\sum_{r\in \cR} m^a_r \delta_{r,j}) = 0, \\
        \pdv{\cL}{\mu_{b(j)}} = 1 - \sum_{t\in\tilde{\cT}}n^b_t = 0, \\
        \pdv{\cL}{\rho_{a(t,w)}} = 1 - \sum_{r\in\tilde{\cR}}m^a_r = 0.
    \end{dcases}
\end{align}

The conditions \eqref{eq:opt_cond_sop} result in the following three conditions:
\begin{align}
    \label{eq:utility_max_ship}
    &(\text{Utility maximization of shipper} ~ b\in\cB_{j})\qquad
    \begin{dcases}
        c^b_t + \lambda^\star_{(t,j)} = \mu^\star_{b(j)} & \text{if $n^{b\star}_t > 0$} \\
        c^b_t + \lambda^\star_{(t,j)} \geq \mu^\star_{b(j)} & \text{if $n^{b\star}_t = 0$}
    \end{dcases}\qquad \forall t\in\tilde{\cT},  \\ %\, \forall b\in\cB_{j},\, \forall j \in \cJ,
    \label{eq:utility_max_drv}
    &(\text{Utility maximization of driver} ~ a\in\cA_{(t,w)})\qquad \begin{dcases}
        c^a_r - \sum_{j\in r}\delta_{r,j}\lambda^\star_{(t,j)} = \rho^\star_{a(t,w)} & \text{if $m^{a\star}_r > 0$} \\
        c^a_r - \sum_{j\in r}\delta_{r,j}\lambda^\star_{(t,j)} \geq \rho^\star_{a(t,w)} & \text{if $m^{a\star}_r = 0$}
    \end{dcases}\qquad \forall r \in \tilde{\cR}, \\ %\, \forall a\in\cA_{(t,w)},\, \forall w \in \cW, \, \forall t\in\cT,
    \label{eq:supply_and_demand}
    &(\text{Market clearing})\qquad
    \begin{dcases}
        \sum_{b\in\cB_{j}}n^{b\star}_t =\sum_{w \in \cW}\sum_{a\in \cA_{(t,w)}}\sum_{r\in \cR} m^{a\star}_r \delta_{r,j} & \text{if $\lambda^\star_{(t,j)} > 0$} \\
        \sum_{b\in\cB_{j}}n^{b\star}_t \leq \sum_{w \in \cW}\sum_{a\in \cA_{(t,w)}}\sum_{r\in \cR} m^{a\star}_r \delta_{r,j} & \text{if $\lambda^\star_{(t,j)} = 0$}
    \end{dcases}\qquad \forall j \in \cJ, \, \forall t \in\cT, 
    % \label{eq:conservation}
    % &(\text{Conservation})\qquad \qquad
    % \sum_{rs \in \cT} y^{rs\star}_a = 1 \qquad \forall a\in\cA.
\end{align}
where $\lambda^\star_{(t,j)}$ can be interpreted as the reward for a task $j$ in the window $t$ at optimum, i.e., $\lambda^\star_{(t,j)} = p^\star_{(t,j)}$. That is, the optimal reward $p^\star_{(t,j)}$ can be obtained by the Lagrangian multiplier with respect to the constraint \eqref{eq:demand_supply}.

The first condition \eqref{eq:utility_max_ship} describes the cost minimization principle of each shipper. In other words, if $b \in \cB_j$ chooses $t$ (i.e., $n^{b\star}_t = 1$), then their disutility equals the minimum disutility $\mu^\star_{b(j)}$. Similarly, the second condition \eqref{eq:utility_max_drv} is the cost minimization principle of each driver, indicating that the disutility of driver $a \in \cA_{(t,w)}$ performing the task bundle $r$ equals the minimum disutility $\rho^\star_{a(t,w)}$ when it is chosen (i.e., $m^{a\star}_r = 1$). 
% This is written in the recursive dynamic programming form, and $\rho^{a\star}_{(k,j)}$ is the minimum cost from states $(k,j)$ to $(K+1,d)$ in the task-chain network.
The third market clearing condition \eqref{eq:supply_and_demand} states that prices are adjusted so as to equilibriate the demand supply of the delivery tasks. Specifically, at market equilibrium, the number of drivers performing task $j$ in $t$ is equal to the number of shippers who request the delivery of task $j$ in $t$ if the reward $\lambda^\star_{(t,j)}$ is larger than zero, and the number of drivers is in excess supply otherwise. %and the manager operates the dedicated vehicle to complete the task

% From the manager's point of view, the manager outsources the tasks to CSD drivers if the reward is cheaper than the operation cost of a dedicated vehicle and performs the tasks by operating dedicated vehicles otherwise. That is, this condition also describes the profit maximization behavior of the manager by choosing between CSD and dedicated drivers.
% Finally, the condition \eqref{eq:conservation} simply states that drivers who express their participation in the CSD system have to be matched with exactly one delivery task. % (i.e., $\overline{c}_{rs} \le w^\star_{rs}$)

That being said, the optimality conditions of the problem [SO] describe the equilibrium under a perfectly competitive market: (i) the maximization of utility for shippers, (i) the maximization of utility for CSD drivers, and (iii) the supply and demand equilibrium (market clearing) condition. As such, the following proposition with the definition of equilibrium in the CSD system holds.
\begin{definition}\label{def:equi}
(Equilibrium matching pattern).
An equilibrium matching pattern under a perfectly competitive market is a tuple of state variables $(\mathbf{m}^\star, \mathbf{n}^\star, \bm{\mu}^\star, \bm{\rho}^\star, \bm{\lambda}^\star)$ that satisfies the conditions \eqref{eq:utility_max_ship},
\eqref{eq:utility_max_drv}, \eqref{eq:supply_and_demand}.
% , and \eqref{eq:conservation}.
\end{definition}
\begin{proposition}\label{prop:1} 
An equilibrium matching pattern realized under the CSD system maximizes the social surplus defined by \eqref{eq:SO_obj}.
\end{proposition}

\section{Fluid-particle decomposition approach}\label{sec:mechanism}
The theoretical analysis in \cref{sec:problem} reveals that a market equilibrium matching pattern achieves social cost minimization. However, achieving the market equilibrium is non-trivial and requires a mechanism design so that shippers/drivers have no incentive to take strategic action. The platform has to solve [SO] for the mechanism design, but it is challenging due to (i) computational complexity and (ii) unobservability of private cost of each shipper/driver. To address these issues, this section presents a hierarchical matching mechanism based on the FPD approach.

%, where in the upper stage we perform task partition to decompose the problem and in the lower stage we solve individual matching problems in a much smaller scale. However, to achieve the market equilibrium and global matching efficiency, we must perform the problem decomposition considering the matching efficiency achieved in the lower stage. To achieve this without explicitly solving the lower-stage sub-problems, we rely on a fluid approximation method that exploits the additive random utility model (ARUM) framework and analytically derives the optimal value functions of the sub-problems. \cite{akamatsu_oyama_csd_2024} originally propose this decomposition method and call it a \textit{fluid-particle decomposition approach}, and we extend it by considering demand/supply elasticities and task bundling.

\subsection{Hierarchical market decomposition}\label{subsec:decompose}
We introduce an aggregate variable $\mathbf{y} = [y_{j,t}]_{j \in \cJ, t \in \tilde{\cT}}$ representing the number of shipping permits for time window $t$ issued for shipper group $\cB_j$, as well as $\mathbf{f} = [f_{(t,w),r}]_{r \in \cR, (t,w) \in \cT\times\cW}$ representing the number of task bundles $r$ allocated to driver group $\cA_{(t,w)}$.
% drivers who travel between OD pair $w$ in $t$ to be matched with bundle $r$. 

% $\mathbf{q} = [q^{(t,j)}_{w}]_{j \in \cJ, w \in \cW, t\in\cT}$ representing the number of tasks $j$ delivered in window $t$ allocated to a set of drivers $\cA_{(t,w)}$ who travel between OD pair $w$ in $t$. Moreover, we introduce $\mathbf{f} = [f_{(t,w),r}]_{r \in \cR, (t,w) \in \cT\times\cW}$ as a variable representing the number of drivers in $\cA_{(t,w)}$ to be matched with bundle $r$.

Then, the problem [SO] can be equivalently formulated as the following integer programming problem:
\begin{align}
    [\text{SO-2}]\qquad
    \min_{\mathbf{m}, \,\mathbf{n}, \,\mathbf{y}, \, \mathbf{f}}\quad&
    \sum_{j\in\cJ} \qty(\sum_{b\in\cB_{j}} \sum_{t\in\tilde{\cT}} n^b_t c^b_t) +
    \sum_{t\in\cT}\sum_{w \in \cW} \qty(
    \sum_{a\in\cA_{(t,w)}}\sum_{r \in \tilde{\cR}}m^a_r c^a_r)
     \label{eq:SO2_obj}\\
    \subto\quad 
    % & \eqref{eq:indv_shipper_consv_all}, \eqref{eq:indiv_driver_consv_all}, \nonumber \\
    &\sum_{t\in\tilde{\cT}} n^b_t = 1 \qquad \forall b \in \cB_{j}, \, \forall j \in \cJ, \tag{\ref{eq:indv_shipper_consv_all}}\\
    &\sum_{r\in \tilde{\cR}} m^a_r = 1 \qquad \forall a \in \cA_{(t,w)}, \, \forall (t,w) \in \cT \times \cW
    \tag{\ref{eq:indiv_driver_consv_all}}\\
    &\sum_{b\in\cB_{j}}n^b_t = y_{j,t} \qquad \forall t \in \tilde{\cT}, \, \forall j\in\cJ\label{eq:ship_eq_permit}\\
    &\sum_{a\in\cA_{(t,w)}}m^a_r = f_{(t,w),r} \qquad \forall r \in \tilde{\cR}, \, \forall (t,w) \in \cW \times \cT,\label{eq:drv_eq_pathflow}\\
    % &\sum_{w \in\cW} q^{(t,j)}_{w} = y_{j,t} \qquad \forall j \in\cJ, \, \forall t \in \cT \label{eq:assign_eq_permit}\\
    % &\sum_{r\in \tilde{\cR}}\delta_{r,j} f_{(t,w),r} \leq q^{(t,j)}_{w} \qquad \forall w \in \cW, \, \forall j\in\cJ,\, \forall t \in \cT,\label{eq:pathflow_leq_assign}\\
    &y_{j,t} \leq \sum_{w \in\cW}\sum_{r\in \tilde{\cR}}\delta_{r,j} f_{(t,w),r} \qquad \forall j \in\cJ, \, \forall t \in \cT,\label{eq:pathflow_leq_assign}\\
    &y_{j,t} \in \Z_{+} \qquad \forall t\in\tilde{\cT}, \, \forall j\in\cJ,\label{eq:y_integer}\\
    % &q^{(t,j)}_{w} \in \Z_{+} \qquad \forall w \in \cW, \,\forall j\in\cJ, \, \forall t\in\cT,\label{eq:q_integer}\\
    &f_{(t,w),r} \in \Z_{+} \qquad \forall r\in\tilde{\cR}, \, \forall (t,w) \in \cW \times \cT\label{eq:f_integer}
\end{align}
where $\Z_{+}$ is the set of all positive integers. 
%$\cA_{w}$ is the set of all drivers whose OD pair of their trips is $(o,d)$ and 
The original supply-demand constraint \eqref{eq:demand_supply} is replaced with constraints \eqref{eq:ship_eq_permit}, \eqref{eq:drv_eq_pathflow} and \eqref{eq:pathflow_leq_assign} using the intermediate variable $\mathbf{y}$ and $\mathbf{f}$. %, $\mathbf{q}$, \eqref{eq:assign_eq_permit}
Note that \eqref{eq:indv_shipper_consv_all} and \eqref{eq:ship_eq_permit}, and \eqref{eq:indiv_driver_consv_all} and \eqref{eq:drv_eq_pathflow} respectively ensure 
% together ensure $N^S_j = \sum_{t \in \tilde{\cT}} y_{j,t}$, and \eqref{eq:indiv_driver_consv_all} and \eqref{eq:drv_eq_pathflow} together ensure $N^D_{(t,w)} = \sum_{r \in \tilde{\cR}} f_{(t,w),r}$.
\begin{align}
    & \sum_{t \in \tilde{\cT}} y_{j,t} = N^S_j \qquad \forall j \in \cJ \label{eq:N_ship_consv}\\
    &\sum_{r \in \tilde{\cR}} f_{(t,w),r} = N^D_{(t,w)} \qquad \forall (t,w) \in \cT \times \cW\label{eq:N_drv_consv}
\end{align}
where $N^S_j$ is the size of shipper group $\cB_j$ and $N^D_{(t,w)}$ is the size of driver group $\cA_{(t,w)}$.

Then, we can hierarchically decompose the problem [SO-2] into a master problem and sub-problems as follows:
\begin{align}
    [\text{SO/Master}]\qquad&\nonumber\\\min_{\mathbf{y}, \mathbf{f}}\quad&
    \sum_{j\in\cJ}z^{S\star}_{j}(\mathbf{y}_j) +
    \sum_{t\in\cT}\sum_{w \in \cW}z^{D\star}_{(t,w)}(\mathbf{f}_{(t,w)})  \label{eq:master_obj}\\ %\mathbf{q},
    \subto\quad& \eqref{eq:pathflow_leq_assign}, \eqref{eq:N_ship_consv}, \eqref{eq:N_drv_consv}\nonumber\\ %\eqref{eq:assign_eq_permit},
    % & \sum_{t \in \tilde{\cT}} y_{j,t} = N^S_j \qquad \forall j \in \cJ \label{eq:N_ship_consv}\\
    % &\sum_{r \in \tilde{\cR}} f_{(t,w),r} = N^D_{(t,w)} \qquad \forall (t,w) \in \cT \times \cW\label{eq:N_drv_consv}\\
    \text{where}\quad&
    z^{S\star}_{j}(\mathbf{y}_j) \equiv \min_{\mathbf{n}_{(j)}} z^S_j(\mathbf{n}_{(j)}|\mathbf{y}_j) \label{eq:optimal_value_func_ship}\\
    & z^{D\star}_{(t,w)}(\mathbf{f}_{(t,w)}) \equiv \min_{\mathbf{m}_{(t,w)}} z^D_{(t,w)}(\mathbf{m}_{(t,w)}|\mathbf{f}_{(t,w)}) \label{eq:optimal_value_func_drv} \\
    [\text{SO/Sub-Shipper$(j)$}]\qquad&\nonumber\\
    \min_{\mathbf{n}_{(j)}} \quad
    & z^S_j(\mathbf{n}_{(j)}|\mathbf{y}_j) \equiv \sum_{b \in \cB_j} \sum_{t\in\tilde{\cT}} n^b_t c^b_t \label{eq:SO_sub_j_obj}\\
    \subto\quad&
    \sum_{t\in\tilde{\cT}} n^b_t = 1 \qquad \forall b\in\cB_j,\label{eq:indiv_shipper_consv_j}\\
    &\sum_{b\in\cB_j}n^b_t = y_{j,t} \qquad \forall t\in\tilde{\cT},\label{eq:y_leq_q_indiv}\\
    &n^b_t\in\{0,1\} \qquad \forall b\in\cB_{j}, \, \forall t\in\tilde{\cT}\label{eq:y_binary_j}. \\
    [\text{SO/Sub-Driver$(t,w)$}]\qquad&\nonumber\\
    \min_{\mathbf{m}_{(t,w)}} \quad
    & z^D_{(t,w)}(\mathbf{m}_{(t,w)}|\mathbf{f}_{(t,w)}) \equiv \sum_{a\in\cA_{(t,w)}}
    \sum_{r \in \tilde{\cR}} m^a_r c^a_r
    \label{eq:SO_sub_od_obj}\\
    \subto\quad&\sum_{r \in \tilde{\cR}} m^a_r = 1
    \qquad \forall a \in \cA_{(t,w)},\label{eq:indiv_driver_consv_od}\\ %\mathbf{q},
    &\sum_{a\in\cA_{(t,w)}}m^a_r = f_{(t,w),r} \qquad \forall r \in \tilde{\cR},\label{eq:m_eq_f_od}\\
    % &\sum_{r \in \tilde{\cR}}\delta_{r,j} f_{(t,w),r} \leq q^{(t,j)}_{w} \qquad \forall j\in\cJ,\label{eq:f_leq_q_od}\\
    &m^a_r \in \{0,1\} \qquad \forall r \in \tilde{\cR}, \, \forall a \in \cA_{(t,w)} \label{eq:ma_binary_Aod}
    % & z^{\star}_{w}(\mathbf{f}) \equiv \max_{\mathbf{m}_{(t,w)}} z_{w}(\mathbf{m}_{(t,w)}|\mathbf{f}) \label{eq:z_star_def}
\end{align}
% Here the constraints \eqref{eq:indiv_driver_consv_od} and \eqref{eq:y_binary_od} are restricted to $\cA_{w}$.
The objective function of the master problem \eqref{eq:master_obj} is defined as the sum of optimal value functions of [SO/Sub-Shipper$(j)$] and [SO/Sub-Driver$(t,w)$], respectively represented by \eqref{eq:SO_sub_j_obj} and \eqref{eq:SO_sub_od_obj}.
%$z^{\star}_{w}(\mathbf{f}) \equiv \max_{\mathbf{m}_{(t,w)}} z_{w}(\mathbf{m}_{(t,w)}|\mathbf{f})$.
% Note that the constraints of the problem [SO/Sub-Driver$(t,w)$] are restricted to $\cA_{w}$, but the differences from that of the problem [SO-2] are slight.

% We now have a master problem and small-scale sub-problems by decomposing the original problem [SO]. However, both of the problems are still nontrivial to solve. Although the sub-problem is a small-scale combinatorial optimization problem, it also involves the private information $c^{b}_{k}$ and $c^a_r$ for both shippers and drivers, which cannot be directly observed by the platform manager. 
% In addition, because the master problem is a task partition problem and simultaneously determines the number of issued shipping permits and that of delivery tasks allocated to each group of drivers, it must be solved before the sub-problems. Nevertheless, the objective function of the master problem \eqref{eq:master_obj} includes the optimal value functions \eqref{eq:optimal_value_func_ship} and \eqref{eq:optimal_value_func_drv} of the sub-problems to achieve global matching efficiency. 
% We will present an extended FPD approach that addresses both of these problems.

\subsection{Fluid approximation of the master problem}
We now have a master problem and small-scale sub-problems by decomposing the original problem [SO].
The master problem of the proposed mechanism determines (i) the number of shipping permits $\mathbf{y}$ to be issued for each shipper group, and (ii) the number of task bundles $\mathbf{f}$ to be allocated to each driver group.
The objective \eqref{eq:master_obj} of the master consists of the optimal value functions of the sub-problems $z^{S\star}$ and $z^{D\star}$, which are obtained as a result of the sub-problems. On the other hand, because the sub-problems are formulated with the quantities $\mathbf{y}$ and $\mathbf{f}$ determined by the master problem as a given, the master problem must be solved before the sub-problems. That is, we in theory must solve the master problem and sub-problems simultaneously.
% In other words, we must evaluate $z^{S\star}_{j}(\mathbf{y})$ and $z^{D\star}_{(t,w)}(\mathbf{f})$ without solving the sub-problems directly.

To address this issue, we extend the fluid-approximation of \cite{akamatsu_oyama_csd_2024}. The core idea of the approach is evaluating $z^{S\star}$ and $z^{D\star}$ as closed-form functions of $\mathbf{y}$ and $\mathbf{f}$ by considering the continuous distributions of shippers' and drivers' utilities. The rationale behind this is that if the auctions described in \Cref{sec:auction} are conducted on a daily basis, the platform manager can collect the information on shippers' and drivers' utilities sufficiently to estimate the distribution, even though the utilities of the shippers and drivers appearing on the day are still unknown. 
Specifically, we rely on an ARUM framework as follows.

\begin{assumption}\label{assumption:utility}
    There are sufficiently large numbers of shippers for each $j \in \cJ$ and drivers for each $(t,w) \in \cT \times \cW$. Then, their perceived costs $c^b_t$ and $c^a_{(i,j)}$ can be approximated by
    \begin{align}
        \label{eq:ARUM_disutility_ship}
        &c^b_t \thickapprox C^S_{j,t} - \zeta^b_{j,t} \qquad &\forall b \in \cB_j, \forall j \in\cJ, \, \forall t\in\tilde{\cT}, \\
        \label{eq:ARUM_disutility_drv}
        &c^a_r \thickapprox C^D_{(t,w),r} - \epsilon^a_{(t,w),r} \qquad &\forall a \in \cA_{(t,w)}, \forall (t,w) \in \cW \times \cT, \, \forall r \in \tilde{\cR}
    \end{align}
    where $C^S_{j,t}$ is the deterministic cost of a shipper shipping their parcel in time window $t\in\cT$ or not using the CSD service ($t=0$), and $C^D_{(t,w),r}$ is the deterministic detour cost for a driver who travels between OD pair $w$ in time window $t$ performing a task bundle $r\in \cR$ or not participating in the system ($r=0$). 
    $\zeta^b_{j,t}$ and $\epsilon^a_{(t,w),r}$ are the unobservable random utilities following joint distributions with finite means that are continuous and independent of $C^S_{j,t}$ and $C^D_{(t,w),r}$, respectively. 
    % The platform manager can observe only $C^S_{j,t}$ and $C^{w}_{(i,j)}$, and 
\end{assumption}
% In particular, we assume the utility distribution to be consistent with an additive random utility (ARUM) model but do not restrict it to a specific model. Therefore, the utility $c^{rs}_{a}$ of driver $a\in\cA_{w}$ performing delivery task $(r, s)\in\cT$ can be expressed as
% \begin{align}
%     c^{rs}_{a} \thickapprox C^{rs}_{w} - \varepsilon^{rs}_{w} \qquad \forall od \in\cW, \, \forall rs \in\cT,
% \end{align}
% where $C^{rs}_{w}$ represents the deterministic detour cost (e.g., free flow travel time) for a driver who travels between OD pair $(o,d)$ and performs a task for $(r,s)$, which the system manager can observe, and $\varepsilon^{rs}_{a}$ represents the unobservable random utility of driver $a$ for performing delivery task $(r,s)$, which follows an arbitrary distribution.

Typically, the deterministic cost $C^D_{(t,w),r}$ of a driver performing a task bundle $r$ is defined as the detour travel time \citep[e.g.,][]{wang2016towards, akamatsu_oyama_csd_2024}, while that cost $C^S_{j,0}$ for a shipper not using the CSD service can represent the cost for requesting a professional delivery service or shipping the parcel themself.
% Let $t_{uv}$ denote the travel time of the shortest path between node pair $(u, v)$, then the deterministic disutilities $C^S_{j,t}$ and $C^D_{(t,w),r}$ are typically defined as:
% \begin{align}
%     \label{eq:Cj0}
%     &C^S_{j,0} \equiv t_{j^+j^-}, \\
%     \label{eq:Cod_ij}
%     &C^D_{(t,w),r} \equiv \sum_{k=1}^{K+1} (t_{j^{-}_{k-1}j^{+}_m} + t_{j^{+}_{m}j^{-}_{m}}) - t_{w},
% \end{align}
% where $o^+ = o^- = o$ and $d^+ = d^- = d$.

Under \Cref{assumption:utility}, the negative surplus functions of shippers and drivers are obtained by the following \textit{expected minimum cost} functions:
\begin{align}
    & \Pi^S_j(\mathbf{v}^S_j) = \mathbb{E}\qty[\min_{t\in\tilde{\cT}} \{v^S_{j,t} - \zeta_{j,t}\}],\label{eq:surplus_ship}\\
    & \Pi^D_{(t,w)}(\mathbf{v}^D_{(t,w)}) = \mathbb{E}\qty[\min_{r\in\tilde{\cR}} \{ v^D_{(t,w),r} - \epsilon_{(t,w),r} \} ], \label{eq:surplus_drv}
\end{align}
where $v^S_{j,t} = C^S_{j,t} + p_{(t,j)}$ and $v^D_{(t,w),r} = C^D_{(t,w),r} - \sum_{j\in \cJ} \delta_{r,j}p_{(t,j)}$ are the generalized deterministic costs. Moreover, for any joint distribution of random utilities $\zeta_{j,t}$ and $\epsilon_{(t,w),r}$, we can find the following \textit{perturbation functions}, or generalized entropy functions \citep{sorensen2022mcfadden, fosgerau2020discrete}, as the convex conjugate of the surplus functions $\Pi^S_j$ and $\Pi^D_{(t,w)}$\citep[see, for example][]{rockafellar1970}:
\begin{align}
    & \cH^S_j(\mathbf{P}^S_j) = \min_{\mathbf{v}^S_{j}} \qty{\mathbf{v}^{S\top}_{j} \mathbf{P}^S_j - \Pi^S_j(\mathbf{v}^S_{j})},\label{eq:entropy_ship}\\
    & \cH^D_{(t,w)}(\mathbf{P}^D_{(t,w)}) = \min_{\mathbf{v}^D_{(t,w)}} \qty{\mathbf{v}^{S\top}_{(t,w)} \mathbf{P}^D_{(t,w)}  - \Pi^D_{(t,w)}(\mathbf{v}^D_{(t,w)})}.\label{eq:entropy_drv}
\end{align}
where $\mathbf{P}^S_j = [P^S_j(t)]_{t\in\tilde{\cT}}$ and $\mathbf{P}^D_{(t,w)}=[P^D_{(t,w)}(r)]_{r\in\tilde{\cR}}$ are probability distributions according to $\zeta_{j,t}$ and $\epsilon_{(t,w),r}$, and they can be interpreted as task partition ratios $P^S_j(t) = y_{j,t}/N^S_j$ and $P^D_{(t,w)}(r) = f_{(t,w),r}/N^D_{(t,w)}$. Thus, we obtain
\begin{align}
    & \hat{\cH}^S_j(\mathbf{y}_j) = N^S_j \cH^S_j(\mathbf{P}^S_j),\label{eq:entropy_ship_y}\\
    & \hat{\cH}^D_{(t,w)}(\mathbf{f}_{(t,w)}) = N^D_{(t,w)} \cH^D_{(t,w)}(\mathbf{P}^D_{(t,w)}).\label{eq:entropy_drv_f}
\end{align}
Then, the following proposition holds.
% Moreover, we treat the numbers $\Theta = (\mathbf{y}, \mathbf{q}, \mathbf{f})$ as continuous variables. 
\begin{proposition}\label{prop:so_a_p}
    Under \Cref{assumption:utility}, the master problem [SO/Master] can be fluidly approximated as 
    \begin{align}
    [\text{SO/Master-A}]\qquad&\nonumber\\
    \min_{\mathbf{y}, \mathbf{f}}\quad& \sum_{j \in \cJ} \qty(\mathbf{C}^{S}_j \cdot \mathbf{y}_j - \hat{\cH}^S_j(\mathbf{y}_j)) + \sum_{t\in\cT}\sum_{w\in \cW} \qty(    \mathbf{C}^{D}_{(t,w)} \cdot \mathbf{f}_{(t,w)} - \hat{\cH}^D_{(t,w)}(\mathbf{f}_{(t,w)}))\label{eq:SOAP_mas_obj}\\
    % \max_{\Pi \in \Delta}\quad& \sum_{j \in \cJ}N^S_j \qty( \mathbf{C}^\top_j \mathbf{P}^S_j - \cH^S_j(\mathbf{P}^S_j)) + \sum_{t\in\cT}\sum_{od\in \cW}N^D_{(t,w)} \qty(\mathbf{C}^\top_{(t,w)} \mathbf{P}^D_{(t,w)} - \cH^D_{(t,w)}(\mathbf{P}^D_{(t,w)}))\label{eq:SOAP_mas_obj}\\
    \subto\quad 
    & \eqref{eq:pathflow_leq_assign},\eqref{eq:N_ship_consv},\eqref{eq:N_drv_consv}\nonumber
    %\eqref{eq:assign_eq_permit},
    % &\sum_{t\in\tilde{\cT}} y_{j,t} = N^S_j \quad \forall j \in \cJ, \label{eq:y_eq_N} \\
    % &\sum_{r\in\tilde{\cR}} f_{(t,w),r} = N^D_{(t,w)} \quad \forall w \in \cW, \, \forall t\in \cT, \label{eq:f_eq_N}
    \end{align}
    % where $N^S_j \equiv |\cB_j|$ and $N^D_{(t,w)} \equiv |\cA_{(t,w)}|$, and we we define the partition ratios $P^S_j(t) = y_{j,t}/N^S_j$ and $P^D_{(t,w)}(r) = f^{r,od}_r/N^D_{(t,w)}$.     
    % $\cH^S_j$ and $\cH^D_{(t,w)}$ are the perturbation functions, which are the convex conjugate of the surplus functions $\Pi^S_j$ and $\Pi^D_{(t,w)}$ of shippers and drivers, respectively:
\end{proposition}
\begin{proof}
    \Cref{app:proof_so_a}.
\end{proof}
The problem [SO/Master-A] is a minimization problem in the perturbed utility form \citep{fudenberg2015stochastic, fosgerau2022perturbed, oyama2022markovian}. As such, the objective function \eqref{eq:SOAP_mas_obj} represents the total cost under the random utility-maximizing behavior of shippers and drivers. 
The perturbed functions $\hat{\cH}^S$ and $\hat{\cH}^D$ capture the heterogeneity in shippers' and drivers' perceived costs and improve the matching efficiency. This is an important difference from previous studies that considers only deterministic costs \citep[e.g.,][]{wang2016towards, soto2017matching}\footnote{In \Cref{sec:experiment}, we benchmark our mechanism against a mechanism considering only the deterministic terms of \eqref{eq:SOAP_mas_obj}, which we call \textit{N0 mechanism}.}.  
% Since the problem's objective is the perturbed utility maximization for shippers and drivers, the partition ratios $\mathbf{P}^S_j$ and $\mathbf{P}^D_{(t,w)}$ are equivalent to the \textit{choice probability distributions} of the ARUMs.
% satisfying
% \begin{align}
%     &\mathbf{P}^S_j \in \arg\min_{\mathbf{P}^S_j} (\mathbf{v}^\top_{j} \mathbf{P}^S_j - H_j(\mathbf{P}^S_j)) \\
%     &\mathbf{P}^D_{(t,w)} \in \arg\min_{\mathbf{P}^D_{(t,w)}} (\mathbf{v}^\top_{(t,w)} \mathbf{P}^D_{(t,w)} - H_{(t,w)}(\mathbf{P}^D_{(t,w)}))
% \end{align}
% Finally, \Cref{prop:so_a_p} holds for any joint distribution of random utilities $\zeta_{j,t}$ and $\epsilon_{(t,w),r}$ \citep{akamatsu_oyama_csd_2024, sorensen2022mcfadden}.

Furthermore, using the convex conjugate duality of $\cH$ and $\Pi$ we can formulate the dual problem of [SO/Master-A].
\begin{proposition}\label{lemma:so_d}
The Lagrangian dual problem of [SO/Master-A] is given as:
\begin{align}
    [\text{SO/Master-D}]\qquad&\nonumber\\
    \max_{\mathbf{p}}\quad&
    G(\mathbf{p}) \equiv \sum_{j \in \cJ} N^S_j \Pi^S_j(\mathbf{p})  +
    \sum_{t\in\cT}\sum_{w \in \cW} N^D_{(t,w)} \Pi^D_{(t,w)}(\mathbf{p}) 
    \label{eq:SOD_obj}
\end{align}
\end{proposition}
\begin{proof}
    \Cref{app:proof_so_d}.
\end{proof}
The dual problem [SO/Master-D] is an unconstrained optimization problem whose objective is the total negative surplus (social cost) of shippers and drivers. The only decision variable is the price vector $\mathbf{p}=[p_{(t,j)}]_{j\in\cJ, t\in\cT}$, which is the dual variable associated with constraint \eqref{eq:pathflow_leq_assign}. % \eqref{eq:assign_eq_permit} and
% The gradient of $G$ with respect to $p_{(t,j)}$ is:
% \begin{align}
%     \nabla_{p_{(t,j)}}G = y^\star_{j,t} - \sum_{w \in \cW}\sum_{r\in\tilde{\cR}} \delta_{r,j} f^{(t,w)\star}_{r}
% \end{align}
% representing the gap between the number of permits issued for task $j$ in $t$ and the total number of tasks for $j$ allocated to drivers who travel in $t$.

\section{Transformation to traffic assignment problem}\label{sec:network}
The fluidly approximated master problem [SO/Master-A] determines continuous variables instead of integer variables, which significantly reduces the computational effort to solve it. However, we still need the enumeration of task bundles $\cR_{(t,w)}$ to perform the task-driver matching. Because the enumeration is prohibitively costly when the number of tasks is large, we herein present an efficient solution algorithm for the master problem, in which we cast the task partition problem as a TAP.

% \subsection{Dual formulation}
% \begin{lemma}\label{lemma:so-d}
%     The Lagrangian dual problem of [SO/Master-A] is obtained as follows:
% \begin{align}
%     [\text{SO/Master-D}]\qquad&\nonumber\\
%     \max_{\mathbf{p}}\quad&
%     \sum_{j \in \cJ} N^S_j \Pi^S_j(\mathbf{p})  +
%     \sum_{t\in\cT}\sum_{w \in \cW} N^D_{(t,w)} \Pi^D_{(t,w)}(\mathbf{p}) 
%     \label{eq:drv_obj}
% \end{align}
% \end{lemma}
% \begin{proof}
%     \Cref{app:proof_so_d}.
% \end{proof}
% Since the dual problem [SO/Master-D] is an unconstrained optimization problem and takes $\mathbf{p}$ as the only decision variable, it is more efficient to solve than the primal problem [SO/Master-A]. However, the computational efficiency depends on the evaluation of the surplus function $\Pi^S_j$ and $\Pi^D_{(t,w)}$ given a price vector $\mathbf{p}$. 
% In particular, the size of the bundle set is $|\cR| = (|\cJ| + 1)^{K+1}$ and generally very large, which imposes the computational burden of the computation of $\Pi^D_{(t,w)}$. We present a dynamic programming approach 

% \subsection{Dynamic programming approach}

\subsection{Matching as network path choice}
\begin{figure}[t]
    \centering
    \includegraphics[width=0.60\textwidth]{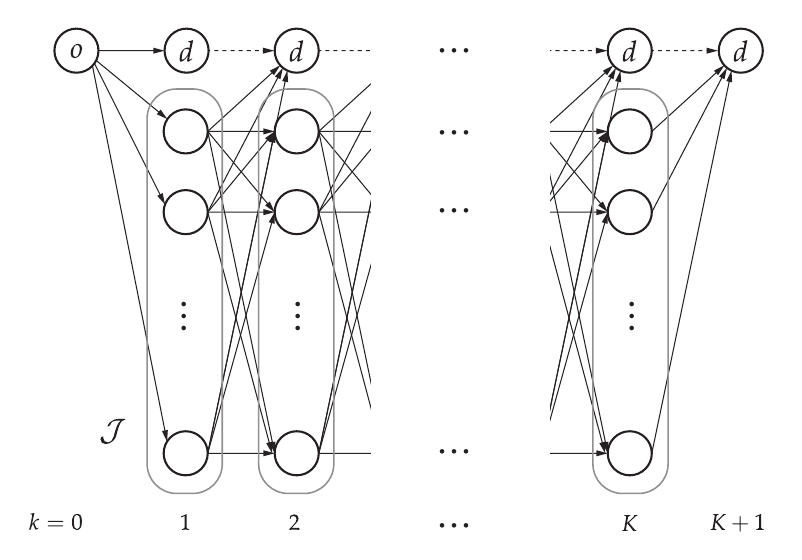}
    \caption{Task-chain network.}
    \label{fig:driver_net}
    \vspace{0.5cm}
    \includegraphics[width=0.95\textwidth]{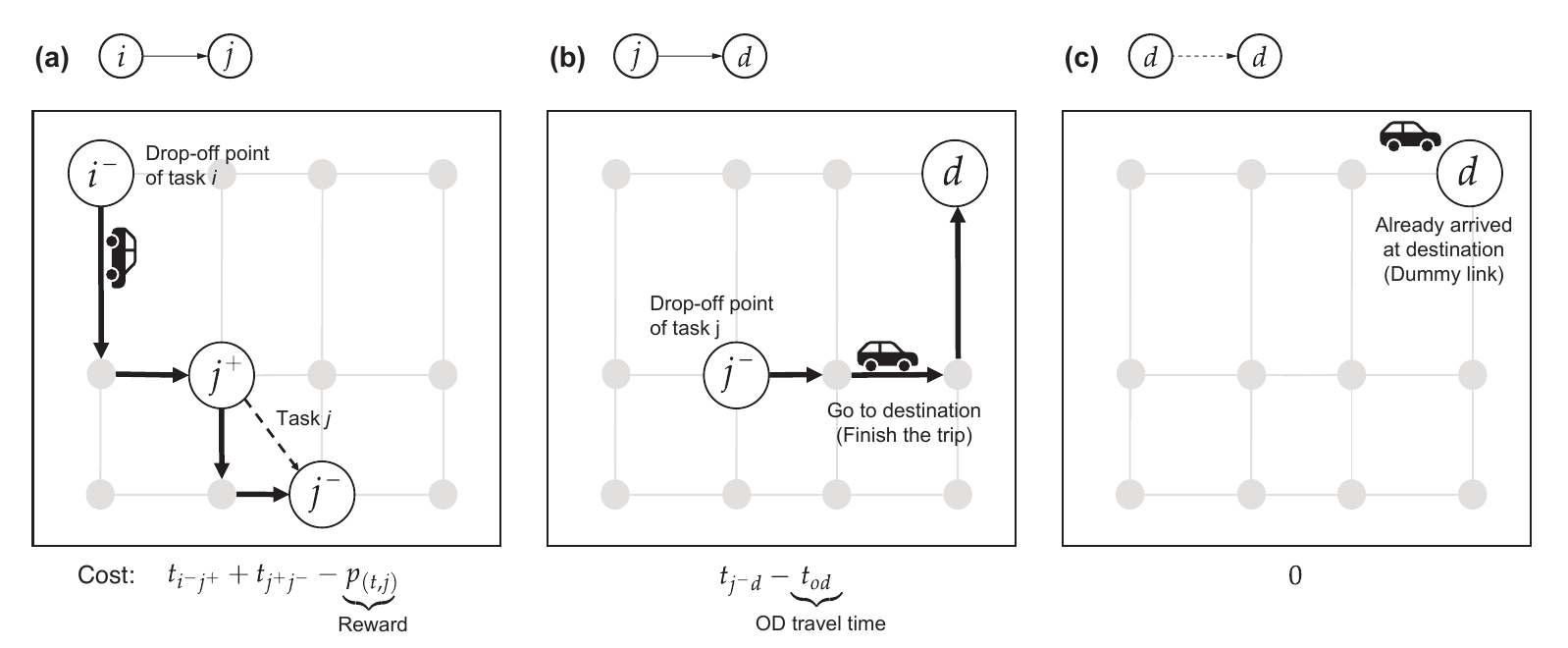}
    \caption{Three types of edges in the network and the associated costs.}
    \label{fig:driver_link}
\end{figure}

% We cast the master problem [SO/Master-A] as a traffic assignment problem on 
We introduce a \textit{task-chain network} of Figure \ref{fig:driver_net}. This network makes $K$ copies of the task set $\tilde{\cJ} \equiv \cJ \cup \{d\}$ and connects each task pair by edge $(i,j) \in \cE = \tilde{\cJ} \times \tilde{\cJ}$. Then, any bundle of tasks can be represented as a path on the network $r = [j_0, j_1, \ldots, j_K, j_{K+1}]$, where $j_0 = o$, $j_{K+1} = d$, and $j_k \in \tilde{\cJ}, \, \forall k \in \{1, \ldots, K\}$. 
When $j_k = d$, the driver finishes their trip after completing $k-1$ tasks, and $j_k = j_{k+1} = \cdots = j_K = d$. The driver's original trip plan without performing a task ($r=\emptyset$) is represented by the path with $j_1 = d$. 
% Then, any task bundle $r \in \tilde{\cR}$ including no-performance ($r=\emptyset$) can be represented as the sequence of $K+1$ edges from $o$ at $k=0$ to $d$ at $k=K+1$, i.e., a path on the network. 
As such, the task bundle-driver matching becomes a path choice problem under this network representation.

However, the number of possible paths $|\tilde{\cR}|$ is about the size of $(|\cJ| + 1)^{K+1}$, which is generally very large and its enumeration is computationally prohibitive. To efficiently model the path choice, we exploit the edge-based approach based on the Markov decision process (MDP). That is, the path choice is represented as a recursive decision process of sequential edge choices for $K+1$ stages. 
To do so, we consider \textit{edge-additive} costs.
% assume that the cost $C^D_{(t,w),r}$ of bundle $r$ is \textit{edge-additive}.

\begin{assumption}\label{assumption:edge_additive}
    The perceived cost of drivers to perform task bundle $r$ is \textit{edge-additive}, that is,
    \begin{align}
        C^D_{(t,w),r} + \epsilon_{(t,w),r} = \sum_{k=0}^{K} \sum_{(i,j) \in \cE} \delta^r_{k,ij} (c^D_{ij} + \epsilon^{(t,w)}_{ij})
    \end{align}
    where $c^D_{ij}$ is the edge cost of $(i,j) \in \cE$ and $\epsilon_{ij}$ is its random term, both of which are independent of $k$, and $\delta^r_{k,ij}$ equals one if path $r$ uses edge $(i,j)$ at $k$, i.e., the $k$-th and $(k+1)$-th tasks of bundle $r$ are $i$ and $j$, and zero otherwise. 
\end{assumption}

\Cref{fig:driver_link} shows the generalized edge costs on the network. For a pair of two tasks $i,j \in \cJ$, the cost is equal to the travel cost $t_{i^-j^+}$ between the drop-off point $i^-$ of task $i$ and the pickup point $j^+$ of task $j$, and $t_{j^+j^-}$ between the pickup and drop-off points $(j^+, j^-)$ of task $j$, minus the reward $p_{(t,j)}$ associated with task $j$. For edge $(j,d)$ connecting $j \in \cJ$ to destination $d$, we subtract the travel cost $t_{od}$ between $(o,d)$ so that the total path cost describes the detour cost for drivers. The network includes dummy edges $(d,d)$ to represent bundles of less than $K$ tasks, and the cost of a dummy edge is zero and deterministic.

Under MDP, the cost-minimizing behavior of drivers is sequentially described at each choice stage $k$. Then, the path choice probability $P^D_{(t,w)}(r)$ is decomposed into 
\begin{align}
   P^D_{(t,w)}(r) = \prod_{k=0}^K \prod_{(i,j) \in \cE} \delta^r_{k,ij} P^D_k(j|i) 
\end{align}
where $P^D_k(j|i)$ is the edge choice probability of $(i,j)$ conditional on state $(k,i)$. In addition, the path-based surplus function $\Pi^D_{(t,w)}$ can be decomposed into $\bm{\pi}_{(t,w)}$, which is recursively formulated for each state $(k,i)$ as 
\begin{align}
    \label{eq:value_func}
    \pi^{(t,w)}_{(k,i)}(\mathbf{p}) = \mathbb{E} \qty[\min_{j\in\tilde{\cJ}} \{ c^D_{ij} - p_{(t,j)} - \epsilon^{(t,w)}_{ij} + \pi^{(t,w)}_{(k+1,j)}(\mathbf{p}) \}  ]. 
    % \quad \forall k \in \{0,\ldots,K+1\}, \, \forall i \in \tilde{\cJ}, \, \forall w \in \cW, \forall t \in \cT,
\end{align}
This is also called the value function of each state and represents the shortest path marginal cost from state $(k,i)$ to the target state $(K+1, d)$ on the network \citep{oyama2019prism, oyama2023capturing}.

\subsection{Markovian traffic assignment}\label{sec:assignment}
Under the network representation, our problem becomes finding of the path flow pattern $\mathbf{f}_{(t,w)}$. 
%defines a state $(k,j)$ as the pair of $k \in \{0,1,\ldots,K,K+1\}$ and $j \in \tilde{\cJ} \equiv \cJ\cup\{d\}$. An edge connecting two states $(k-1,i)$ and $(k,j)$ represents the driver performing the $m$-th task $j$ after completing task $i$. Therefore, a task bundle $r = [j_0, j_1, \ldots, j_K, j_{K+1}]$ can be seen as a path on this network from states $(0,o)$ to $(K+1,d)$, and $f_{(t,w),r}$ as a path flow. 
% However, the number of possible paths $|\cR|$ is $(|\cJ| + 1)^{K+1}$ and generally very large. To efficiently compute the flow on the network, we exploit a link-based Markovian traffic assignment framework \citep{akamatsu1996cyclic, oyama2022markovian}. 
Moreover, as we exploit a edge-based traffic assignment framework, the decision variable is reduced to the edge flow pattern $\mathbf{x}_{(t,w)}$ that satisfies: % \citep{akamatsu1996cyclic, oyama2022markovian}
\begin{align}
    \label{eq:link_path_flow}
    x^{(t,w)}_{k,ij} = \sum_{r \in \tilde{\cR}} \delta^r_{k,ij} f_{(t,w),r}.
\end{align}
At each state $(k,i)$, the conditional edge choice probability $P^D_k(j|i)$ partitions the sum of upstream flows into edge flows:
\begin{align}
    x^{(t,w)}_{k,ij} = P^D_k(j|i) \sum_{h\in\tilde{\cJ}} x^{(t,w)}_{k,hi},
\end{align}
and in particular, at the source state $(0,o)$,
\begin{align}
    x^{(t,w)}_{0,oj} = P_0(j|o) N^D_{(t,w)}.
\end{align}

With this vector of edge flows and under the edge-additive assumption of random utilities, the perturbation function $\hat{\cH}^D_{(t,w)}(\mathbf{f}_{(t,w)})$ is also redefined as
\begin{align}
    \hat{\cH}^D_{(t,w)}(\mathbf{f}_{(t,w)}) = \bar{\cH}^D_{(t,w)}(\mathbf{x}_{(t,w)}) = \sum_{k=0}^K \sum_{i\in\cJ} \bar{\cH}^{(t,w)}_{k,i}(\mathbf{x}^{(t,w)}_{k,i})
\end{align}
where $\mathbf{x}^{(t,w)}_{k,i} = [x^{(t,w)}_{k,ij}]_{\forall j: (i,j)\in\cE}$.

% We then define the perturbed function $\cH^D_{(t,w)}$ as the sum of perturbations across edges:
% \begin{align}
%     \cH^D_{(t,w)}(\mathbf{P}^D_{(t,w)}) = \sum_{k=0}^{K} \sum_{(i,j) \in \cE} \delta^r_{k,ij} H_{k,ij}(P^D_k(j|i))
% \end{align}

% \subsection{Traffic assignment problem formulation}
Using this network representation, the problem [SO/Master-A] can be equivalently formulated as the following traffic assignment problem (TAP):
\begin{align}
    [\text{SO/Master-TAP}]\qquad&\nonumber\\
    \max_{\mathbf{y}, \mathbf{x}}\quad&
    \sum_{j \in \cJ} \qty( \sum_{t \in \tilde{\cT}} C^S_{j,t} y_{j,t} - \hat{\cH}^S_j(\mathbf{y}_j)) +
    \sum_{t\in\cT}\sum_{w \in \cW} \qty(\sum_{k=0}^{K}\sum_{i,j\in \tilde{\cJ}} c^D_{ij} x^{(t,w)}_{k,ij}  - \bar{\cH}^D_{(t,w)}(\mathbf{x}_{(t,w)}))
    \label{eq:TAP_obj}\\
    \subto\quad
    & \sum_{t \in \tilde{\cT}} y_{j,t} = N^S_j \qquad \forall j \in \cJ \tag{\ref{eq:N_ship_consv}}\\
    % &  \eqref{eq:N_ship_consv},\nonumber\\
    %\eqref{eq:assign_eq_permit},
    %&\mathbf{A}\mathbf{x}_{w} = \mathbf{q}_{w} \\
    &\sum_{i\in\tilde{\cJ}} x^{(t,w)}_{k-1,ij} - \sum_{l\in\tilde{\cJ}} x^{(t,w)}_{k,jl} = \eta_{k,j} N^D_{(t,w)}  \quad \forall k \in \{0,\ldots,K+1\}, \, \forall j \in \cJ, \label{eq:x_flow_consv}\\
    % & \sum_{k=0}^{K}\sum_{i\in\tilde{\cJ}} x^{(t,w)}_{k,ij} \leq q^{(t,j)}_{w} \quad \forall j\in\cJ, \, \forall w \in \cW, \, \forall t \in \cT \label{eq:x_flow_cap}
    & y_{j,t} \leq \sum_{w\in\cW}\sum_{k=0}^{K}\sum_{i\in\tilde{\cJ}} x^{(t,w)}_{k,ij}  \quad \forall j\in\cJ, \, \forall t \in \cT \label{eq:x_flow_cap} %\, \forall w \in \cW
\end{align}
where constraint \eqref{eq:x_flow_consv} is the flow conservation law on the network, and
$\eta = [\eta_{k,j}]$ is a two-dimensional vector whose elements equal zero except for $\eta_{K+1,d} = 1$ and $\eta_{0,o} = -1$. This replaces constraint \eqref{eq:N_drv_consv} of [SO/Master-A].
% Constraint \eqref{eq:x_flow_consv} is the flow conservation law on the network, and 
Constraint \eqref{eq:x_flow_cap} is the demand-supply constraint ensuring that the number of shipping permits $y_{j,t}$ of time window $t$ for shippers with a task $j$ is less than the number of drivers who perform tasks for $j$ in $t$, which replaces \eqref{eq:pathflow_leq_assign} of [SO/Master-A].
% These replace [SO/Master-A]'s constraints \eqref{eq:N_drv_consv} and \eqref{eq:pathflow_leq_assign}, respectively. 
The problem [SO/Master-TAP] is a (capacitated) Markovian, or link-based stochastic, traffic assignment problem \citep{akamatsu1996cyclic, akamatsu1997decomposition, Baillon2008MCA, oyama2019prism,oyama2022markovian}. More generally, it can be viewed as a perturbed utility traffic assignment problem \citep{fosgerau2022perturbed, yao2023perturbed}.

% \subsection{Dual formulation}
Furthermore, the dual problem of [SO/Master-TAP] is obtained by re-formulating [SO/Master-D] as:
\begin{align}
    [\text{SO/Master-TAP-D}]\qquad&\nonumber\\
    \max_{\mathbf{p}}\quad&
    G(\mathbf{p}) \equiv \sum_{j \in \cJ} N^S_j \Pi^S_j(\mathbf{p})  +
    \sum_{t\in\cT}\sum_{w \in \cW} N^D_{(t,w)} \pi^{(t,w)}_{(0,o)}(\mathbf{p}) 
    \label{eq:dual_obj_edge}
    % \text{where}\quad& 
    % \pi^{(t,w)}_{(k,i)}(\mathbf{p}) = \mathbb{E} \qty[\min_{j\in\tilde{\cJ}} \{ c^D_{ij} - p_{(t,j)} - \epsilon^{(t,w)}_{ij} + \pi^{(t,w)}_{(k+1,j)}(\mathbf{p}) \}  ], \quad \forall k \in \{0,\ldots,K+1\}, \, \forall i \in \tilde{\cJ}, \, \forall w \in \cW, \forall t \in \cT,
    % \label{eq:drv_exp_min_cost}
\end{align}
% \begin{proof}
%     \Cref{app:proof_so_d} (or \Cref{app:proof_so_d_edge}).
% \end{proof}
wherein the driver surplus function $\Pi^D_{(t,w)}$ in the objective \eqref{eq:SOD_obj} of [SO/Master-D] is replaced by $\pi^{(t,w)}_{(0,o)}$. As defined in \eqref{eq:value_func}, the value function $\pi^{(t,w)}_{(k,j)}$ of each state is recursively formulated, and at the source state, $\pi^{(t,w)}_{(0,o)}$ implicitly represents the expected minimum cost among all possible paths on the network, i.e., all task bundles, thus being equivalent to $\Pi^D_{(t,w)}$.
% With the edge-based formulation, the surplus (or the expected minimum cost) function $\pi^{(t,w)}_{(k,j)}$ is defined for each state in the network in the recursive form. At the source state, $\pi^{(t,w)}_{(0,o)}$ implicitly represents the expected minimum cost of all possible paths, i.e., all task bundles, thus being equivalent to $\Pi^D_{(t,w)}$. 
This recursively-formulated $\pi^{(t,w)}_{(k,j)}$ can be solved through a dynamic programming approach for any ARUMs \citep{oyama2022markovian}.

\subsection{Examples}
\subsubsection{Example 1: The multinomial logit model}
We herein show as an example the specific case of the multinomial logit (MNL) model for the fluid approximation. To do so, we assume the unobserved costs $\zeta_{j,t}$ and $\epsilon^{(t,w)}_{ij}$ of shippers and drivers follow type I extreme value (EV) distributions: $\zeta_{j,t} \stackrel{iid}{\sim} \text{EV}(0, \theta)$ and $\epsilon^{(t,w)}_{ij} \stackrel{iid}{\sim} \text{EV}(0, \phi)$.

For the MNL models, it is well-known that the perturbation functions are defined as Shannon's entropy functions: 
% For the shipper model, $\cH^S_j(\mathbf{y}_j)$ becomes
% \begin{align}
%     \cH^S_j(\mathbf{y}_j) = -\frac{1}{\theta} \sum_{t\in\tilde{\cT}} y_{j,t} \ln  \frac{y_{j,t}}{N^S_j}, %= N^S_j \cH^S_j(\mathbf{P}^S_j)
% \end{align}
% and accordingly, the surplus function $\Pi^S_j(\mathbf{p})$ is the following log-sum function:
% \begin{align}
%     \Pi^S_j(\mathbf{p}) = -\frac{1}{\theta} \ln \sum_{t\in\tilde{\cT}} \exp{-\theta (C^S_{j,t} + p_{(t,j)})} \label{eq:ship_surplus_logit} 
% \end{align}
\begin{align}
    &\hat{\cH}^S_j(\mathbf{y}_j) = -\frac{1}{\theta} \sum_{t\in\tilde{\cT}} y_{j,t} \ln  \frac{y_{j,t}}{N^S_j} \\ %= N^S_j \cH^S_j(\mathbf{P}^S_j)
    &\bar{\cH}^D_{(t,w)}(\mathbf{x}_{(t,w)}) = \sum_{k=0}^K \sum_{i\in\cJ} \hat{\cH}^{(t,w)}_{k,i}(\mathbf{x}^{(t,w)}_{k,i})
     = -\frac{1}{\phi} \sum_{k=0}^{K} \sum_{i\in\cJ} \sum_{j\in\cJ} x^{(t,w)}_{k,ij}\ln  \frac{x^{(t,w)}_{k,ij}}{\sum_{h\in\cJ}x^{(t,w)}_{k-1,hi}} %= N^D_{(t,w)} \cH^D_{(t,w)}(\mathbf{P}^D_{(t,w)})
\end{align}
and accordingly, the surplus functions are the following log-sum functions:
\begin{align}
    &\Pi^S_j(\mathbf{p}) = -\frac{1}{\theta} \ln \sum_{t\in\tilde{\cT}} \exp{-\theta (C^S_{j,t} + p_{(t,j)})} \label{eq:ship_surplus_logit} \\
    &\pi^{(t,w)}_{(k,i)}(\mathbf{p}) = -\frac{1}{\phi} \ln \sum_{j\in\tilde{\cJ}} \exp{-\phi (c^D_{ij} - p_{(t,j)} + \pi^{(t,w)}_{(k+1,j)}(\mathbf{p}))} \label{eq:drv_surplus_logit}
\end{align}
% These are closed-form solutions, and 
The recurrence relation of \eqref{eq:drv_surplus_logit} reduces to the system of linear equations \citep{akamatsu1996cyclic, Fosgerau2013RL}, and the solution can be simply computed as shown in \Cref{sec:loading}.

Moreover, the demand (network loading) functions with the MNL probabilities are
\begin{align}
    &y^\star_{j,t} = N^S_j \exp{-\theta (C^S_{j,t} + p_{(t,j)} - \Pi^S_j(\mathbf{p}))} \label{eq:demand_y}\\
    % &x^{(t,w)\star}_{0,oj} = N^D_{(t,w)} \exp{-\phi (c^D_{ij} - p_{(t,j)} + \pi^{(t,w)}_{(1,j)}(\mathbf{p}) - \pi^{(t,w)}_{(0,o)}(\mathbf{p}))}\\
    &x^{(t,w)\star}_{k,ij} = \qty(\sum_{h\in\cJ}x^{(t,w)\star}_{k,hi}) \cdot \exp{-\phi (c^D_{ij} - p_{(t,j)} + \pi^{(t,w)}_{(k+1,j)}(\mathbf{p}) - \pi^{(t,w)}_{(k,i)}(\mathbf{p}))}, \label{eq:demand_x}
\end{align}
with which we can compute the task partition patterns given the price pattern $\mathbf{p}$.

\subsubsection{Example 2: The nested logit model}
We then show another example of the nested logit (NL) model that captures the correlation among unobserved utilities.
Since both the models of shippers and drivers involve \textit{out-out} options, we can consider the nest structure focusing on participation in CSD. That is, we assume that a shipper (driver) first chooses among two nests, i.e., whether or not to participate in the CSD system, and if participating, then the shipper (driver) chooses a time window for delivery (task bundle to perform). The participation nest for the shipper model includes all time window alternatives in $\cT$, and that for the driver model consists of all task bundles in $\cR$ for the driver model.
% Note that different specifications of the unobserved cost distributions lead to other models. One can specify a nested logit model for the shippers' choice, as the choice involves (i) whether or not to outsource a task and (ii) which time window to outsource, thus the correlation in unobserved utilities of the time window alternatives.

The NL formulations for the shipper model are:
% For the NL model, the perturbation functions are defined as
\begin{align}
    &\hat{\cH}^S_j(\mathbf{y}_j) = -\frac{1}{\theta} \qty(y_{j,0} \ln \frac{y_{j,0}}{N^S_j} + \sum_{t\in\cT} y_{j,t} \ln \frac{\sum_{t\in\cT} y_{j,t}}{N^S_j}) - \frac{1}{\hat{\theta}}    
    \sum_{t\in\cT} y_{j,t} \ln  \frac{y_{j,t}}{\sum_{t\in\cT} y_{j,t}} \label{eq:ship_entropy_NL}\\ %= N^S_j \cH^S_j(\mathbf{P}^S_j)
    &\Pi^S_j(\mathbf{p}) = -\frac{1}{\theta} \ln \qty{\exp(-\theta C^S_{j,0}) + \exp(-\theta\Pi^S_{j,\cT}(\mathbf{p}))}
    \label{eq:ship_surplus_NL} \\
    &\Pi^S_{j,\cT}(\mathbf{p}) = -\frac{1}{\hat{\theta}}\ln \sum_{t\in\cT} \exp{-\theta (C^S_{j,t} + p_{(t,j)})} \label{eq:ship_surplus_in_NL}
\end{align}
where $\Pi^S_{j,\cT}(\mathbf{p})$ represents the expected minimum cost given $\mathbf{p}$ among shipping time windows. Two scales $\theta$ and $\hat{\theta}$ are associated with nest choice and time window choice within the participation nest, respectively.
As such, the demand functions are
\begin{align}
    &y^\star_{j,0} = N^S_j \exp{-\theta (C^S_{j,0} - \Pi^S_j(\mathbf{p}))} \label{eq:y_out_NL}\\
    % &y^\star_{j,t} = N^S_j \exp{-\hat{\theta} (C^S_{j,t} + p_{(t,j)}) - (\theta - \hat{\theta}) \Pi^S_{j,\cT}(\mathbf{p}) - \theta \cS_{j}(\mathbf{p})} \label{eq:demand_y}\\
    &y^\star_{j,\cT} = N^S_j \exp{-\theta (\Pi^S_{j,\cT}(\mathbf{p}) - \cS_{j}(\mathbf{p}))} \label{eq:y_in_NL}\\
    &y^\star_{j,t} = y^\star_{j,\cT} \exp{-\hat{\theta} (C^S_{j,t} + p_{(t,j)} - \Pi^S_{j,\cT}(\mathbf{p}))} \label{eq:y_NL}
\end{align}
where $y^\star_{j,\cT}$ is the total number of shippers who outsource the delivery of task $j$, while $y^\star_{j,0}$ is the number of shippers who do not participate in the CSD system.

The NL formulations for the driver model are obtained similarly. For clarity, we first show the path-based formulations as follows:
\begin{align}
    &\hat{\cH}^D_{(t,w)}(\mathbf{f}_{(t,w)}) = -\frac{1}{\phi} \qty(f^{(t,w)}_{\emptyset} \ln \frac{f^{(t,w)}_{\emptyset}}{N^D_{(t,w)}} + \sum_{r\in\cR} f_{(t,w),r} \ln \frac{\sum_{t\in\cR} f_{(t,w),r}}{N^D_{(t,w)}}) - \frac{1}{\hat{\phi}}    
    \sum_{r\in\cR} f_{(t,w),r} \ln  \frac{f_{(t,w),r}}{\sum_{r\in\cR} f_{(t,w),r}} \label{eq:drv_entropy_NL} \\ 
    &\Pi^D_{(t,w)}(\mathbf{p}) = -\frac{1}{\phi} \ln \qty{\exp(-\phi C^D_{(t,w),\emptyset}) + \exp(-\phi\Pi^D_{(t,w),\cR}(\mathbf{p}))}
    \label{eq:drv_surplus_NL} \\
    &\Pi^D_{(t,w),\cR}(\mathbf{p}) = -\frac{1}{\hat{\phi}}\ln \sum_{r\in\cR} \exp{-\hat{\phi} (C^D_{(t,w),r} - p_{t,r})} \label{eq:drv_surplus_in_NL}
\end{align}
where $\Pi^D_{(t,w),\cR}(\mathbf{p})$ represents the expected minimum cost given $\mathbf{p}$ among all possible task bundles. Like the shipper model, the NL driver model has two scales $\phi$ and $\hat{\phi}$ associated with nest and bundle choices, respectively. The demand functions are
\begin{align}
    &f^{(t,w)\star}_{\emptyset} = N^D_{(t,w)} \exp{-\phi (C^D_{(t,w),\emptyset} - \Pi^D_{(t,w)}(\mathbf{p}))} \label{eq:f_out_NL}\\
    &f^{(t,w)\star}_{\cR} = N^D_{(t,w)} \exp{-\phi (\Pi^D_{(t,w),\cR}(\mathbf{p}) - \Pi^D_{(t,w)}(\mathbf{p}))} \label{eq:f_in_NL}\\
    &f^{(t,w)\star}_r = f^{(t,w)\star}_{\cR} \exp{-\hat{\phi} (C^D_{(t,w),r} - p_{t,r} - \Pi^D_{(t,w),\cR}(\mathbf{p}))} \label{eq:f_NL}
\end{align}
On the task-chain network, $f^{(t,w)}_{\emptyset} = x^{(t,w)}_{0,od}$, $f^{(t,w)}_{\cR} = N^D_{(t,w)} - x^{(t,w)}_{0,od}$ and $C^D_{(t,w),\emptyset} = c_{od}$. Simply by removing edge $(o,d)$ from the network, the MNL path choice model with scale $\hat{\phi}$ represents the within-nest choice model of the NL model, thereby
\begin{align}
    % &\cH^D_{(t,w)}(\mathbf{x}_{(t,w)}) = - \frac{1}{\phi} \qty(x^{(t,w)}_{0,od} + (N^D_{(t,w)}-x^{(t,w)}_{0,od}) \ln (N^D_{(t,w)}-x^{(t,w)}_{0,od})) - \frac{1}{\hat{\phi}} \sum_{k=0}^{K} \sum_{i\in\cJ} \sum_{j\in\cJ} x^{(t,w)}_{k,ij}\ln  \frac{x^{(t,w)}_{k,ij}}{\sum_{h\in\cJ}x^{(t,w)}_{k-1,hi}}\\
    &\pi^{(t,w)}_{(k,i)}(\mathbf{p}) = -\frac{1}{\hat{\phi}} \ln \sum_{j\in\tilde{\cJ}} \exp{-\hat{\phi} (c^D_{ij} - p_{(t,j)} + \pi^{(t,w)}_{(k+1,j)}(\mathbf{p}))} \label{eq:drv_surplus_rec_NL}\\
    % &x^{(t,w)\star}_{k,ij} = \qty(\sum_{h\in\cJ}x^{(t,w)\star}_{k,hi}) \cdot \exp{-\theta_{D,i} (c^D_{ij} - p_{(t,j)} + \pi^{(t,w)}_{(k+1,j)}(\mathbf{p}) - \pi^{(t,w)}_{(k,i)}(\mathbf{p}))}, \label{eq:demand_x}
    &x^{(t,w)\star}_{k,ij} = \qty(\sum_{h\in\cJ}x^{(t,w)\star}_{k,hi}) \cdot \exp{-\hat{\phi} (c^D_{ij} - p_{(t,j)} + \pi^{(t,w)}_{(k+1,j)}(\mathbf{p}) - \pi^{(t,w)}_{(k,i)}(\mathbf{p}))}, \label{eq:x_NL}
\end{align}
and $\Pi^D_{(t,w),\cR} = \pi^{(t,w)}_{(0,o)}$ and $f^{(t,w)\star}_{\cR} = \sum_{j\in\cJ}x^{(t,w)\star}_{0,oj}$ hold. 
As such, the surplus and demand functions for the NL model can also be decomposed into edge-based formulations on the network. 

\section{Solution algorithm of the master problem}\label{sec:algorithm}
\subsection{Accelerated gradient descent method}
% By casting the master problem as a traffic assignment problem, we reduce the computational order of the problem. The original [SO] involves $(|\cB| \times (|\cT| + 1)) + (|\cA| \times (|\cJ| + 1)^{K+1})$ decision variables, which is reduced to $(|\cB| \times (|\cT| + 1)) + (|\cA| \times (K+1) \times (|\cJ| + 1)^{2})$. 
We exploit the dual formulation of [SO/Master-TAP-D] to develop an efficient solution algorithm. The primal problem [SO/Master-TAP] is a stochastic TAP but has capacity constraints. Solving a capacitated traffic assignment problem (CTAP) has been challenging, and solution algorithms applicable to large-scale networks have not been developed in the literature \citep{nie2004models, Meng2008}.
% and requires the iterative process of network loading.  

In contrast, the dual problem [SO/Master-TAP-D] is unconstrained with the price pattern $\mathbf{p}$ being the only decision variable, which allows us to apply existing fast non-linear optimization algorithms. This study presents an efficient dual algorithm based on the accelerated gradient descent (AGD) method \citep{Nesterov1983, oyama2022markovian}. The AGD method is a first-order (gradient) method and requires less memory space compared to second-order methods. It also achieves the known complexity bound $\Omega(1/k^2)$ , that is, in the worst case, any iterative method based solely on the function and gradient evaluations cannot achieve a better accuracy than $\Omega(1/k^2)$ at iteration $k$. Although dual type AGD methods are presented in the context of traffic equilibrium assignment \citep{oyama2022markovian, yao2023perturbed}, this study is the first to propose to develop an AGD-based method for CTAP.

The gradient of $G$, the objective function \eqref{eq:dual_obj_edge} of [SO/Master-TAP-D], is obtained as follows:
\begin{align}
    \nabla_{p_{(t,j)}}G = y^\star_{j,t} - \sum_{w \in \cW}\sum_{k=0}^{K}\sum_{i\in\cJ} x^{(t,w)\star}_{k,ij}
\end{align}
This represents the gap between the number of permits for time window $t$ to be issued for task $j$ and the number of tasks for $j$ to be partitioned for drivers traveling in time window $t$. Therefore, the gradients have an intuitive interpretation as they are equal to \textit{excess demand} for the delivery of tasks in the time window.

We can obtain these flows $\mathbf{y}^\star$ and $\mathbf{x}^\star$ by performing network loading \citep{akamatsu1996cyclic, oyama2022markovian}, and we write these as:
\begin{align}
    &\mathbf{y}^\star_j = \cF^S_j(\mathbf{p}) \\
    &\mathbf{x}^\star_{(t,w)} = \cF^D_{(t,w)}(\mathbf{p})
\end{align}
where $\cF$ is a demand (network loading) function: see for examples \eqref{eq:demand_y} and \eqref{eq:demand_x} for the MNL model case and \eqref{eq:y_out_NL}-\eqref{eq:y_NL} and \eqref{eq:x_NL} for the NL model case. 
% The network loading algorithm is a sub-module of our algorithm and is described in \cref{sec:loading}.

Then, we apply the fast iterative shrinkage-thresholding algorithm (FISTA) that integrates the proximal operator into the AGD method and improves the convergence rate to $\cO(1/k^2)$ where $k$ is the iteration number \citep{beck2009fast}. The solution algorithm for our problem [SO/Master-TAP-D] is summarized as follows.\\

\noindent [Accelerated Gradient Descent (AGD)]
\renewcommand{\theenumi}{\arabic{enumi}}
\renewcommand{\labelenumi}{\textit{Step \theenumi}:}
\renewcommand{\theenumii}{\arabic{enumii}}
\renewcommand{\labelenumii}{\textit{Step \theenumi-\theenumii}:}
\begin{enumerate}
    \setlength{\leftskip}{0.5cm}
    \item \textbf{Initialization}: Initialize price vector $\mathbf{p}^{(0)}$ and step size $\gamma$, and set iteration counter $m := 0$, adaptive restart counter $l := 0$, and auxiliary price $\hat{\mathbf{p}}^{(0)} := \mathbf{p}^{(0)}$. 
    % input step size update rate $\nu$, and
    % and minimum iteration before restart $\underline{l}$ % momentum acceleration variable $\alpha_0 := 0$
    \item \textbf{Network loading}: Compute $\mathbf{y}^{(m)}_j$ for all $j \in \cJ$ and $\mathbf{x}^{(m)}_{(t,w)}$ for all $(t,w) \in \cT \times \cW$:
    \begin{align}
        &\mathbf{y}^{(m)}_j = \cF^S_j(\mathbf{p}^{(m)}) \nonumber\\
        &\mathbf{x}^{(m)}_{(t,w)} = \cF^D_{(t,w)}(\mathbf{p}^{(m)})\nonumber
    \end{align}
    % \item \textbf{Backtracking}: Find the smallest non-negative integer $\iota$ such that, with $\gamma_{m+1} = \nu^{\iota} \gamma_m$,
    % \begin{align}
    %     & G(\mathbf{p}^{(m)} + \gamma_{m+1} \nabla_{p}G(\mathbf{p}^{(m)})) \geq G(\mathbf{p}^{(m)}) + \frac{\gamma_{m+1}}{2}||\nabla_{p}G(\mathbf{p}^{(m)})||^2
    % \end{align}
    \item \textbf{Updating}: Update the current solution using a momentum:
    \begin{align}
        &p^{(m+1)}_{(t,j)} = \hat{p}^{(m)}_{(t,j)} + \gamma \nabla_{p_{(t,j)}}G(\mathbf{p}^{(m)}) = \hat{p}^{(m)}_{(t,j)} + \gamma \qty(y^{(m)}_{j,t} - \sum_{w \in \cW}\sum_{k=0}^{K}\sum_{i\in\cJ} x^{(t,w)(m)}_{k,ij}) \\
        & \tau_{l+1} = \frac{1 + \sqrt{1 + 4\tau^2_l}}{2} \\
        & \hat{\mathbf{p}}^{(m+1)} = \mathbf{p}^{(m+1)} + \frac{\tau_l - 1}{\tau_{l+1}} (\mathbf{p}^{(m+1)} - \mathbf{p}^{(m)})
    \end{align}
    \item \textbf{Adaptive restart}: If the restart criterion holds, then $l := 0$. Otherwise, $l := l + 1$.
    \item \textbf{Convergence test}: If the convergence criterion holds, stop the algorithm. Otherwise, set $m := m + 1$ and return to \textit{Step 1}.
\end{enumerate}

The efficiency of the AGD method relies on network loading (i.e., the calculation of gradients) in Step 2. Therefore, the next subsection explains the efficient network loading algorithm.

% This significant reduction in the number of unknown variables can make the solution more stable.

% The dual problem [SO/Master-TAP-D] has two major advantages compared to the primal problem. Firstly, while the decision variables of the primal [SO/Master-TAP] are $\mathbf{y} \in \mathbb{R}^{|\cJ| \times (|\cT| + 1)}_+$, $\mathbf{q} \in \mathbb{R}^{|\cJ| \times |\cT| \times |\cW|}$, and $\mathbf{x} \in \mathbb{R}^{|\cT| \times |\cW| \times (K+1) \times (|\cJ| + 1)^{2}}$, that of the dual [SO/Master-TAP-D] is only $\mathbf{p} \in \mathbb{R}^{|\cJ| \times |\cT|}$. This means significant reductions in the number of decision variables and thus the required memory capacity for computation.

% By the optimality conditions, the demand functions for the loading are:
% \begin{align}
%     &y^\star_{j,t} = (\cH^\prime_{j,t})^{-1}(C^S_{j,t} + p_{(t,j)} - \mu_j) \\
%     &x^{(t,w)\star}_{k,ij} = (\cH^{(t,w)\prime}_{k,ij})^{-1}(c^D_{ij} + p_{(t,j)} + \rho^{(t,w)}_{k+1,j} - \rho^{(t,w)}_{k,i})
% \end{align}
% Furthermore, the AGP method is combined with backtracking to adjust step size  \citep{Beck2009fista} and with adaptive restart schemes to modify momentum near the optimal solution \citep{Donoghue2015adaptive}.

\subsection{Network loading algorithm}\label{sec:loading}
The computation of the shipper flow pattern $\mathbf{y}$ requires only a single-stage calculation of shippers' time-window choice probabilities $\mathbf{P}^S$ given a price pattern $\mathbf{p}$. Then, multiplying this probability $P^S_j(t)$ by the number of shippers $N^S_j$ yields $y_{j,t}$ (see \eqref{eq:demand_y} and \eqref{eq:y_out_NL}-\eqref{eq:y_NL} for the MNL and NL models, respectively). The driver model is not as simple as the shipper model, because the choice set (i.e., the set of task bundles) for drivers cannot be easily enumerated. Therefore, we apply the Markovian traffic assignment (MTA) algorithm \citep{Dial1971, akamatsu1996cyclic, oyama2019prism, oyama2022markovian} that efficiently computes network traffic flows without enumerating paths. The MTA algorithm consists of two main phases: (i) the backward computation of the surplus (expected maximum utility) functions and (ii) the forward computation of edge flows, as described as follows (we focus on a specific driver group $(t,w)$ and omit the subscript for simplicity). \\

\noindent [Markovian traffic assignment (MTA)]
\begin{enumerate}
    \setlength{\leftskip}{0.5cm}
    \item \textbf{Initialization}: Input a reward pattern $\mathbf{p}$ and the number of drivers $N_{(t,w)}$. Define the exponential of surplus functions $\mathbf{V} = [\exp(-\phi \pi_{(k,j)})]$ and edge flows $\mathbf{x} = [x_{k,ij}]$, and initialize their entries to all zero.
    % $V_{(k,j)} = \exp(\pi^D_{(k,j)}) = 0$ for all $k \in \{0, 1, \ldots, K+1\}$ and all $j \in \cJ$, and link flows $x_{k,ij}$ for all 
    \item \textbf{Backward computation}: 
    \begin{enumerate}
        \setlength{\leftskip}{0.5cm}
        \item Update $V_{(k,d)} = 1$ for all $k \in \{0, 1, \ldots, K+1\}$, and then set $k = K$.
        \item For all $j \in \cJ$, compute:
        \begin{align}
            V_{(k,j)} = \sum_{j' \in \tilde{\cJ}} \exp{-\phi (c^D_{jj'} - p_{(t,j')})} V_{(k+1, j')}
        \end{align}
        \item If $k = 0$, go to Step 3. Otherwise, update $k := k-1$ and go back to \textit{Step 2-2}.
    \end{enumerate}
    \item \textbf{Forward computation}: 
    \begin{enumerate}
        \setlength{\leftskip}{0.5cm}
        \item Set $k = 1$, and compute for all $j \in \tilde{\cJ}$:
        \begin{align}
            x_{0,oj} = N_{(t,w)} \exp{-\phi (c^D_{oj} - p_{(t,j)})} V_{(1, j)} / V_{(0, o)}.
        \end{align}
        \item For all $ij \in \cE$, compute:
        \begin{align}
            x_{k,ij} = \sum_{h \in \tilde{\cJ}} x_{k-1, hi} \exp{-\phi (c^D_{ij} - p_{(t,j)})} V_{(k+1, j)} / V_{(k, i)}
        \end{align} 
        \item If $k=K$, finish the algorithm. Otherwise, update $k := k+1$ and go back to \textit{Step 3-2}.
    \end{enumerate}
\end{enumerate}

Note that although the algorithm presented above is for the MNL model, the NL model requires only an additional step in \textit{Step 3-1} to first partition the flow $N_{(t,w)}$ into participating/non-participating drivers.
This MTA algorithm on the task-chain network is very efficient because both forward and backward computation can be performed by vector multiplication and only requires $K$ times computation. This algorithm is a classical method for link-based stochastic network loading in transportation science \citep{Dial1971, akamatsu1996cyclic, oyama2019prism}, but the novelty of this paper lies in effectively utilizing this algorithm in solving the CSD matching problem.

\section{Solution algorithm of the sub-problems}\label{sec:auction}
% \section{Auction mechanisms for particle matching}\label{sec:auction}
The VCG market mechanism \citep{vickrey1961counterspeculation, clarke1971multipart, groves1973incentives} is known as a mechanism that satisfies the two desirable properties of truth-telling and efficiency \citep{cramton2006, vohra2011mechanism}, and we exploit it to solve our sub-problems.

\subsection{Shipper market}
First, we design the VCG mechanism for shipper markets as follows.\\

\noindent [VCG/Shipper]
\begin{enumerate}
    \setlength{\leftskip}{0.5cm}
    \item Each shipper $b \in \cB_j$ submits a bid $s^b_t$ for each time window $t$, indicating their willingness to pay for the delivery in the window. %their value of the window.
    \item The platform assigns a shipping permit to each driver to maximize the social surplus with respect to the declared values, by solving [SO/Sub-Shipper($j$)] with $c^b_t$ in \eqref{eq:SO_sub_j_obj} replaced with $-s^b_t$: %, that is, 
    \begin{align}
        \max_{\mathbf{n}_{(j)}} \quad
        & z^S_j(\mathbf{n}_{(j)}|\mathbf{s},\mathbf{y}) \equiv \sum_{b \in \cB_j} \sum_{t\in\tilde{\cT}} n^b_t s^b_t
        \label{eq:SO_sub_j_obj_bid}\\
        \subto\quad&
        \sum_{t\in\tilde{\cT}} n^b_t = 1 \qquad \forall b\in\cB_j,\tag{\ref{eq:indiv_shipper_consv_j}}\\
        &\sum_{b\in\cB_j}n^b_t = y_{j,t} \qquad \forall t\in\tilde{\cT},\tag{\ref{eq:y_leq_q_indiv}}\\
        &n^b_t\in\{0,1\} \qquad \forall b\in\cB_{j}, \, \forall t\in\tilde{\cT}\tag{\ref{eq:y_binary_j}}
        % &\eqref{eq:indiv_shipper_consv_j},\, \eqref{eq:y_leq_q_indiv},\, \eqref{eq:y_binary_j}, \, \eqref{eq:n_binary}. \nonumber
    \end{align}
    \item The price $p_b$ that shipper $b \in \cB_j$ pays is determined by
    % \begin{align}
    %     \label{eq:vcg_reward_ship}
    %     p_b(\mathbf{s}_b, \mathbf{s}_{-b}) = 
    %     \max_{\mathbf{n}_{-b}} z^S_{j,-b}(\mathbf{n}_{-b}|\mathbf{s}_{-b},\mathbf{y}) - 
    %     \qty(z^{S\star}_j(\mathbf{s}_b, \mathbf{s}_{-b}) - \sum_{t\in\tilde{\cT}} n^{b\star}_t s^b_t )
    % \end{align}
    \begin{align}
        \label{eq:vcg_reward_ship}
        p_b(\mathbf{s}_b, \mathbf{s}_{-b}) = 
        \max_{\mathbf{n}_{-b}} z^S_{j,-b}(\mathbf{n}_{-b}|\mathbf{s}_{-b},\mathbf{y}) - 
        \sum_{b'\in \cB_j\setminus\{b\}}\sum_{t\in\tilde{\cT}} n^{b'\star}_t(\mathbf{s}) s^{b'}_t
    \end{align}
    where %$z^{S\star}_j(\mathbf{s}) \equiv  \max_{\mathbf{n}_{(j)}} z^S_j(\mathbf{n}_{(j)}|\mathbf{s},\mathbf{y})$ and %the optimal value of the objective and
    $\mathbf{n}^{\star}(\mathbf{s}) = [n^{b\star}_{t}(\mathbf{s})]_{t\in\tilde{\cT}, b\in\cB_j}$ is the optimal matching pattern under bid $\mathbf{s} = (\mathbf{s}_b, \mathbf{s}_{-b})$ for the problem in Step 2, and the subscript ``$-b$'' represents all shippers in $\cB_j$ except for $b$. %, and $\mathbf{s}_{-b}$ represents the vector of bids excluding for $b$.
    The first and second terms of the RHS of \eqref{eq:vcg_reward_ship} are the surplus of shippers $-b$ achieved without and with the participation of $b$, respectively.
    % when the matching problem in Step 2 is performed without driver $a$. 
\end{enumerate}
Herein, Eq. \eqref{eq:vcg_reward_ship} implies that each shipper will pay based on the loss in social surplus caused by the shipper's participation, i.e., the compensation to the opportunity loss of someone else. Under this payment, the VCG mechanism holds efficiency properties.

\begin{proposition}\label{prop:VCG_ship}
    The VCG mechanism for the shipper markets  [VCG/Shipper] is truthful and efficient.
\end{proposition}
\begin{proof}
    \Cref{sec:proof_VCG_ship}.
\end{proof}

This property allows us to collect the information on privately perceived costs of shippers by performing the auctions in the sub-problems.

\subsection{Driver market}
% Scalable edge-based auction mechanism
Using the task-chain network representation, we reformulate the sub-problem [SO/Sub-Driver$(t,w)$]:
\begin{align}
    [\text{SO/Sub-Driver-Edge$(t,w)$}]&\nonumber\\
    \min_{\mathbf{m}_{(t,w)}} \quad
    & z^D_{(t,w)}(\mathbf{m}_{(t,w)}|\mathbf{x}_{(t,w)}) \equiv \sum_{a\in\cA_{(t,w)}}
    \sum_{k=0}^{K} \sum_{ij\in\cE}
    m^a_{k,ij} c^D_{ij}
    \label{eq:SO_sub_od_edge_obj}\\
    \subto\quad
    &\sum_{i\in\tilde{\cJ}} m^a_{k-1,ij} - \sum_{l\in\tilde{\cJ}} m^a_{k,jl} = \eta_{k,j}  \quad \forall k \in \{0,\ldots,K+1\}, \, \forall j \in \cJ,\label{eq:m_edge_flow_consv}\\
    & \sum_{a\in\cA_{(t,w)}} m^a_{k,ij} = x^{(t,w)}_{k,ij}  \quad \forall k \in \{0,\ldots,K+1\}, \, \forall ij \in \cE,
    \label{eq:m_edge_xcap}\\
    & m^a_{k,ij} \in \{0,1\}  \quad \forall k \in \{0,\ldots,K+1\}, \, \forall ij \in \cE, \, \forall a \in \cA_{(t,w)},\label{eq:m_edge_binary}
    % &  \sum_{a\in\cA_{(t,w)}} \sum_{k=0}^K \sum_{i\in\cJ} m^a_{k,ij} \leq q^j_{(t,w)} \label{eq:m_edge_qcap}
\end{align}
where $\mathbf{m}_{(t,w)} = [m^a_{k,ij}]$ is a binary variable that takes $1$ when driver $a \in \cA_{(t,w)}$ is matched with edge $ij$ for $k$-th choice and 0 otherwise, and $\mathbf{x}_{(t,w)}$ is the task partition pattern that is determined in the master TAP problem.
% Note that \eqref{eq:m_edge_qcap} is always satisfied when \eqref{eq:m_edge_xcap} is satisfied due to the capacity constraint \eqref{eq:x_flow_cap} of [SO/Master-TAP]:
% \begin{align}
%     \sum_{a\in\cA_{(t,w)}} \sum_{k=0}^K \sum_{i\in\cJ} m^a_{k,ij} \leq \sum_{k=0}^K \sum_{i\in\cJ} x^{(t,w)}_{k,ij} \leq q^j_{(t,w)}.
% \end{align}
% As such, we can ignore \eqref{eq:m_edge_qcap} when solving [SO/Sub-Driver($t,w$)]. This reduces the problem [SO/Sub-Driver($t,w$)] to a multi-agent shortest path problem with edge capacity constraints.

The platform opens an auction that distributes delivery tasks among each driver group $\cA_{(t,w)}$. In this study, a driver may perform a bundle of delivery tasks, but the number of possible bundles is generally large. Therefore, we utilize the above edge-based formulation to design a scalable VCG auction mechanism.\\
% We introduce a scalable edge-based auction mechanism, which we call a dynamic programming (DP) auction.
% The DP auction consists of sequential $K$-stages auctions. In each stage, drivers submit bids to the edges connected from their states. 
% \ldots

\noindent [VCG/Driver]
\begin{enumerate}
    \setlength{\leftskip}{0.5cm}
    \item Each CSD driver $a \in \cA_{(t,w)}$ submits a bid $\tilde{c}^a_{ij}$ for task pair $(i,j) \in \cE$, indicating their willingness to accept for performing the task $j$ given their location $i^-$. % (i.e., perceived edge costs in the task-chain network)
    \item The platform matches delivery tasks to each driver to minimize the social cost with respect to the declared costs $\tilde{\mathbf{c}}$, by solving [SO/Sub-Edge($t,w$)] with $c^a_{ij}$ in \eqref{eq:SO_sub_od_edge_obj} replaced with $\tilde{c}^a_{ij}$: %, that is, 
    \begin{align}
        \min_{\mathbf{m}_{(t,w)}}\quad&
        z^D(\mathbf{m}_{(t,w)}|\tilde{\mathbf{c}},\mathbf{x}) \equiv \sum_{a\in\cA_{(t,w)}}
        \sum_{k=0}^{K} \sum_{ij\in\cE}
        m^a_{k,ij} \tilde{c}^a_{ij} \label{eq:SO_sub_edge_b}\\
        \subto\quad
        &\eqref{eq:m_edge_flow_consv}, \eqref{eq:m_edge_xcap}, \eqref{eq:m_edge_binary}.
    \end{align}
    \item The reward $w_a$ that CSD driver $a\in\cA_{(t,w)}$ receives is determined by
    % \begin{align}
    %     \label{eq:vcg_reward_drv_edge}
    %     w_{a}(\mathbf{b}_a, \mathbf{b}_{-a}) = 
    %     \min_{\mathbf{m}_{-a}} z^D_{-a}(\mathbf{m}_{-a}|\mathbf{b}_{-a},\mathbf{x}) - 
    %     \qty(z^{D\star}(\mathbf{b}_a, \mathbf{b}_{-a}) - \sum_{k=0}^{K} \sum_{ij\in\cE}
    %     m^{a\star}_{k,ij} b^a_{ij})
    % \end{align}
    \begin{align}
        \label{eq:vcg_reward_drv_edge}
        p_{a}(\tilde{\mathbf{c}}_a, \tilde{\mathbf{c}}_{-a}) = 
        \min_{\mathbf{m}_{-a}} z^D_{-a}(\mathbf{m}_{-a}|\tilde{\mathbf{c}}_{-a},\mathbf{x}) 
        - \sum_{a'\in\cA_{(t,w)}\setminus\{a\}}\sum_{k=0}^{K} \sum_{ij\in\cE}
        m^{a'\star}_{k,ij}(\tilde{\mathbf{c}}) \tilde{c}^{a'}_{ij}
    \end{align}
    where %$z^{D\star}(\mathbf{b}) \equiv  \min_{\mathbf{m}_{(t,w)}} z^D(\mathbf{m}_{(t,w)}|\mathbf{b},\mathbf{x})$ and 
    $\mathbf{m}^{\star}(\tilde{\mathbf{c}}) = [m^{a\star}_{k,ij}(\tilde{\mathbf{c}})]_{i,j\in\tilde{\cJ},k\in\cK,a\in\cA_{(t,w)}}$ is the optimal matching pattern under bid pattern $\tilde{\mathbf{c}}$.
    The subscript ``$-a$'' denotes all drivers except for $a$, and thus $\mathbf{b}_{-a}$ represents the vector of bids excluding $a$, and $\min_{\mathbf{m}_{-a}} z^D_{-a}(\cdot)$ is the social cost value achieved without participation of $a$. 
    % Therefore, \eqref{eq:vcg_reward_drv_edge} states that the reward is equivalent to the contribution of $a$ to reducing the social cost.
\end{enumerate}
Eq. \eqref{eq:vcg_reward_drv_edge} implies that each driver is rewarded based on the driver's contribution to the reduction of social cost, and the reward would be zero if the driver did not contribute to the system. 
%This is an important difference in rewarding to the literature because most existing CSD matching systems deterministically compensate drivers according to their detour costs \citep{archetti2016vehicle, wang2016towards} or decide compensation to maximize platform profit \citep{ccinar2024role}.
The edge-based VCG mechanism for drivers' markets also holds efficiency properties.
\begin{proposition}\label{prop:VCG_drv}
    The VCG mechanism for the driver markets [VCG/Driver] is truthful and efficient.
\end{proposition}
\begin{proof}
    \Cref{sec:proof_VCG_drv}.
\end{proof}

Again, this property allows us to collect the information on privately perceived costs of drivers by performing auctions in the market.

\section{Reduction of computational complexity}
Our approach drastically reduces the computational complexity of the matching problem. Before presenting numerical results, we discuss how much our approach theoretically reduces the complexity and required memory capacity compared to the original matching problem. We refer to \cref{fig:flow} to see how our approach transformed the problem and reduced its size.

The original matching problem [SO] simultaneously decides matchings between drivers and task bundles and between shippers and delivery time windows. The decision variables are $\mathbf{m}$ and $\mathbf{n}$ with sizes of $|\cA| \times |\tilde{\cR}|$ and $|\cB| \times |\tilde{\cT}|$, which can be tremendously large in urban areas due to the growth of $|\cA|$ and $|\cB|$. The number of task bundles $|\tilde{\cR}|$ exponentially increases as $|\cJ|$ and $K$ increase. 

To make this problem tractable, we first apply the FPD approach that decomposes [SO] into the fluidly approximated master problem [SO/Master-A] and small-scale sub-problems [SO/Sub-Shipper($j$)] and [SO/Sub-Driver($t,w$)], by grouping individual shippers by task pickup-delivery location pairs ($\cB \rightarrow \{\cB_j\}_{j\in\cJ}$) and individual drivers by OD pairs and time windows ($\cA \rightarrow \{\cA_{(t,w)}\}_{(t,w)\in\cT \times \cW}$). The decision variables of the master problem are then $\mathbf{y}$ and $\mathbf{f}$ with sizes of $|\cJ| \times |\tilde{\cT}|$ and $|\cT| \times |\cW| \times |\tilde{\cR}|$. 

Then, we transform the master problem into a TAP on the task-chain network [SO/Master-TAP]. This transformation exploits Markovian traffic assignment formulation to address the computational complexity of task-bundling. That is, we do not need to enumerate $\tilde{\cR}$ but implicitly calculate $\mathbf{f}$ based on edge flows $\mathbf{x}$. Furthermore, we propose solving its dual problem [SO/Master-TAP-D] whose decision variable is just $\mathbf{p}$ with a size of $|\cT| \times |\cJ|$. The AGD method finds $\mathbf{p}$ that maximizes the (fluidly approximated) social surplus while satisfying the demand-supply constraint. 

The sub-problems are much smaller problems than the original problem [SO]. In addition, all the sub-problems are completely independent of each other and can be solved in parallel. The driver-side sub-problem [SO/Sub-Driver($t,w$)] is transformed using the task-chain network representation into [SO/Sub-Edge($t,w$)] and can be solved without bundle generation. It is worth noting that the network size for [SO/Sub-Edge($t,w$)] depends on the result of the task partition of the master problem. If only a few tasks category are assigned to a driver group $|\cA_{(t,w)}|$, only a small sub-network of \cref{fig:driver_net} need to be considered and thus the decision variable of [SO/Sub-Edge($t,w$)] becomes very small.

To summarize, our approach transformed the original intractable large-scale problem [SO] into [SO/Master-TAP-D], only finding $\mathbf{p}$ by an efficient AGD method, and independent small-scale sub-problems [SO/Sub-Shipper($j$)] and [SO/Sub-Edge($t,w$)]. By doing so, we drastically reduce computational complexity and save memory space.

Importantly, our approach does not only scale down the size of the problem. The master problem evaluates the optimal value functions of the sub-problems as its objective, by exploiting the RUM framework, so that global matching efficiency is achieved. As such, our approach allows to obtain an accurate solution of the large-scale CSD matching problem with high speed.

\section{Numerical experiments}\label{sec:experiment}
This section presents numerical experiments to demonstrate the performance of the proposed mechanism and algorithm, in terms of both computational efficiency and approximation accuracy. 

We consider the original solution of [SO] as the ground-truth benchmark, where the perceived costs $c^a_r$ and $c^b_t$ of all drivers and shippers are perfectly known. That is, this solution achieves globally optimal matching. We refer to this mechanism as \textit{a naive mechanism with perfect information} (\textbf{N1}).

In addition, to discuss the impact of the observations of drivers' and shippers' private utilities through auctions, we consider another benchmark mechanism. Here, we assume that information on drivers' utilities is unavailable, and the platform solves matching by solving [SO] using only the deterministic cost $C^D_{(t,w),r}$ and $C^S_{j,t}$. We call this mechanism \textit{a naive mechanism without private information} (\textbf{N0}).

Both N1 and N0 mechanisms are solved as a linear programming (LP) model. We use an LP solver of Gurobi Optimizer version 10.0.1 with the default parameter settings and record CPU time in seconds.

Then, we compare the proposed FPD mechanism with the benchmarks. %We use the MNL model for the experiments, that is, 
The FPD mechanism solves shipper-driver matching by utilizing the distributions of shippers' and drivers' private utilities, which are assumed to have been estimated using information collected during past auctions. As described in the paper, the master task-partition problem is solved by the AGD method described in \cref{sec:algorithm}, and the sub-problems are solved by the same LP solver as the benchmarks. We use the MNL model for this experiment, and terminate the AGD method when the iteration reaches 1000 times or when the following three criteria are all satisfied: the maximum gradient $G$ is below $0.1$; the maximum relative difference between $\mathbf{p}^{(m)}$ and $\mathbf{p}^{(m+1)}$ is below $10^{-4}$; and the relative objective difference is below $10^{-6}$.

We perform the experiments in a Python 3.9 environment on a machine with a 12-core Apple M2 Pro and 32 GB of memory, and implement the algorithms by writing our own code. 

\subsection{Performance indicators}
We evaluate the performance of the FPD mechanism by using the following indicators:

\begin{itemize}
    \item \textbf{CPU time}. The CPU time of the FPD mechanism is the sum of the runtime of the AGD method to solve [SO/Master-TAP-D], the average computing time taken to solve [SO/Sub-Shipper($j$)] by the LP solver, and that for [SO/Sub-Driver($t,w$)]. The reasoning behind taking the average over the sub-problems is that the sub-problems are completely independent of each other and can be solved in parallel.
    The CPU times of the N1 and N0 mechanisms are the runtimes to solve [SO] by the LP solver. 
    Note that since the path enumeration is computationally prohibitive, we solve [SO] and [SO/Sub-Driver($t,w$)] using the edge-based formulation on the task-chain network (see \Cref{app:so_edge})\footnote{This means that the naive mechanisms also benefit from the edge-based formulation on the task-chain network. We show that even with this reformulation of [SO], the naive mechanisms still require large computational costs, and our approach drastically reduces the costs.}.
    \item \textbf{Social surplus}. The optimal objective value $z^\star_1$ of [SO] achieved by the N1 mechanism is the ground-truth social surplus that we aim to accurately approximate. 
    The N0 mechanism also solves [SO] but without private information. We define the objective achieved by N0 as $\bar{z}^\star_0$, and measure its error compared to the ground truth by the true relative difference: $(z^\star_1 - z^\star_0)/z^\star_1$.     
    The social surplus achieved by the FPD mechanism is the sum of the optimal objective values of all sub-problems, i.e., $z^\star_{\text{FPD}} = \sum_{j \in \cJ}  z^{S\star}_{j} + \sum_{t \in \cT} \sum_{w \in \cW} z^{D\star}_{(t,w)}$. We measure the approximation error by $(z^\star_1 - z^\star_{\text{FPD}})/z^\star_1$.
    \item \textbf{Price/reward patterns}.     
    We additionally compare the price/reward patterns between mechanisms. The optimal reward pattern of the N1 mechanism $\mathbf{p}^\star_1$ is obtained as the Lagrangian multipliers by solving [SO]. We obtain $\mathbf{p}^\star_1$ for the N0 mechanism in the same way but without private information. The reward pattern of the FPD mechanism $\mathbf{p}^\star_{\text{FPD}}$ is obtained directly by solving the master problem as described in \cref{sec:algorithm}. We measure the approximation error of the reward pattern by the mean of relative differences $\mathbb{E}[(\mathbf{p}^\star_1 - \mathbf{p}^\star_{\text{FPD}})/\mathbf{p}^\star_1]$ and $\mathbb{E}[(\mathbf{p}^\star_1 - \mathbf{p}^\star_0)/\mathbf{p}^\star_1]$.
\end{itemize}

\subsection{Data generation procedure}
The experiments use Winnipeg network data of \citet{trans}, which contains 1,052 nodes, 2,836 links, 147 zones, and 4,345 possible zone pairs for drivers' OD pairs or tasks' pickup-delivery locations.
Based on this network, we generate datasets for the experiments by the following procedure:
\renewcommand{\theenumi}{\arabic{enumi}}
\renewcommand{\labelenumi}{\theenumi.}
\begin{enumerate}
    \item Randomly generate $|\cW|$ different driver OD pairs and $|\cJ|$ different task OD pairs from the Winnipeg network zone pairs.
    \item Randomly decide the OD pair and time window of each of all $|\cA|$ drivers and the pickup-delivery location pair for each of all $|\cB|$ shippers, so that there is at least one driver for every $(t,w) \in \cT \times \cW$ and one shipper for every $j \in \cJ$.
    \item Compute the true perceived costs:  $c^b_t = C^S_{j,t} - \zeta^b_t$ for all shippers $b \in \cB$ and time windows $t \in \cT$; $c^a_{ij} = c^D_{ij} + \epsilon^a_{ij}$ for all drivers $a \in \cA_{(t,w)}$ and edges $ij \in \cE$. As this experiment assumes the MNL model, we randomly draw them from i.i.d. type I EV distributions  $\zeta^b_t \sim \text{EV}(0, \theta)$, and $\epsilon^a_{ij} \sim \text{EV}(0, \phi)$.
\end{enumerate}
We generate 20 datasets for each set of parameters tested in the experiments.
The true private utility values generated in \textit{Step 3} are directly used in the N1 mechanism to obtain the socially optimal solution. We assume that the FPD mechanism only knows the distributional parameters $\theta$ and $\phi$, and the N0 mechanism does not have any information on private utility, i.e., $\zeta^b_t = \epsilon^a_{ij} = 0$. This assumption of N0 is common with the previous matching models that consider deterministic compensation and neglect the stochastic behavior of shippers/drivers \citep[e.g.,][]{wang2016towards, soto2017matching}.

\subsection{Results}
Defining the default setting as $(|\cA|, |\cB|, |\cT|, |\cW|, |\cJ|, K, \theta, \phi) = (5000, 5000, 4, 10, 10, 2, 1.0, 1.0)$, we see how the CPU time and approximation errors change depending on parameter values. We mainly focus on the effects of (1) the numbers of drivers and shippers ($|\cA|,  |\cB|$), (2) the number of task pickup-delivery location pairs ($|\cJ|$), and (3) scale parameter values ($\theta$, $\phi$).

\begin{figure}[t]
    \centering
    \includegraphics[width=\textwidth]{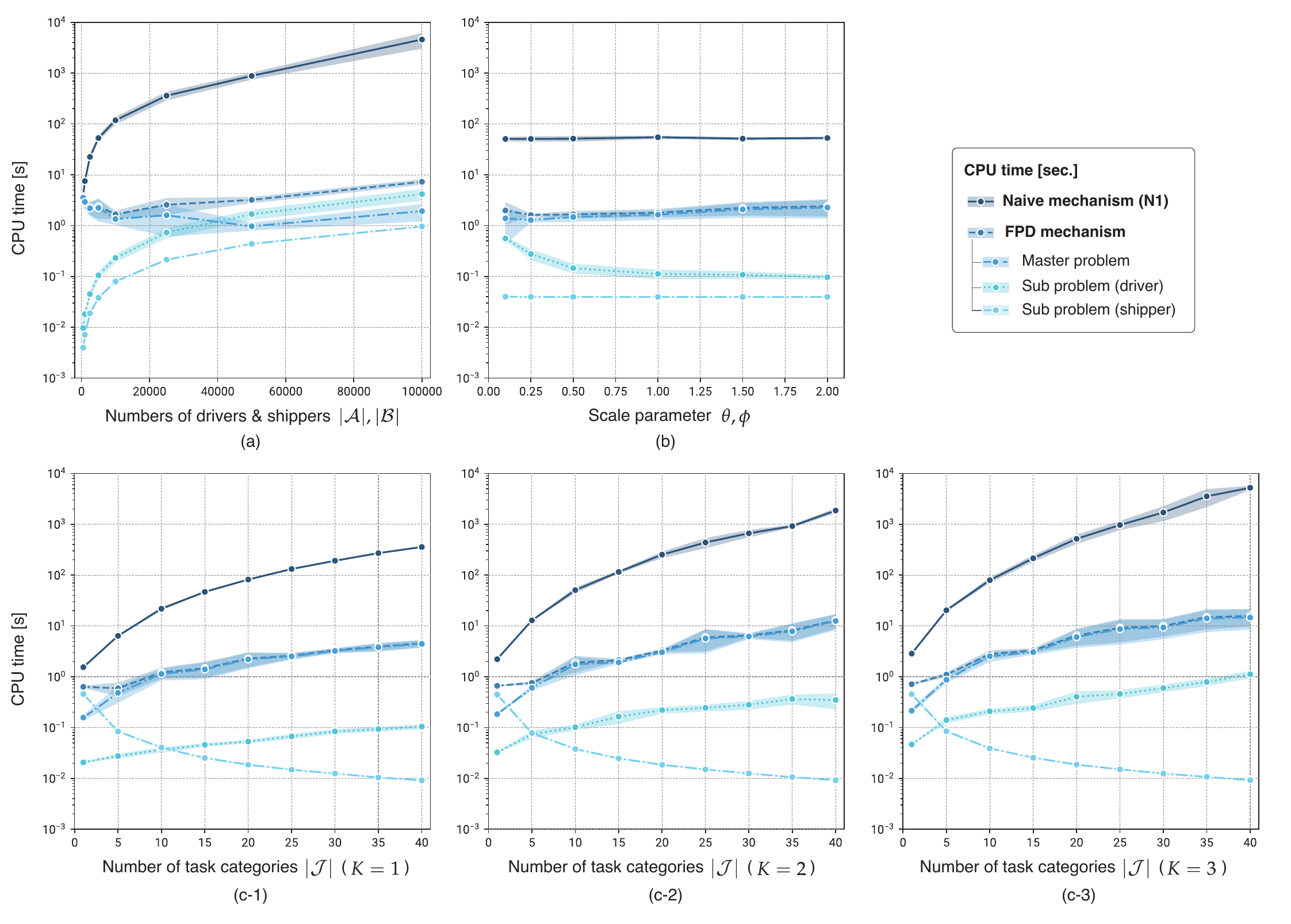}
    \caption{CPU time.}
    \label{fig:cpu_time}
\end{figure}

\paragraph{CPU time.} \cref{fig:cpu_time} shows the CPU times required for N1 and the FPD mechanisms, with averages and standard deviations. First, the increase in the numbers of drivers $|\cA|$ and shippers $|\cB|$ directly enlarges the problem size and requires more CPU time for the N1 mechanism. In contrast, the CPU time of the FPD mechanism is reduced with $|\cA|$ and $|\cB|$. Because our master problem is cast as a flow-independent TAP (as in \cref{sec:assignment}), its computational efficiency does not depend on $|\cA|$ and $|\cB|$. Then the increase in $|\cA|$ and $|\cB|$ tends to lead to a faster convergence of the AGD method. The sizes of the sub-problems increase, but their effects are not large enough to impact the total CPU time. 
Even at $|\cA| = |\cB| = 100000$, the FPD mechanism requires only 7.3 seconds on average, which is about 700 times faster than the N1 mechanism with 5008.7 seconds required on average. This result clarifies a substantial reduction of required computational effort by the proposed approach.

Second, as shown in panel (b) of \cref{fig:cpu_time}, the effect of the scale parameter ($\theta$, $\phi$) on CPU times is not large. The convergence of the master problem is faster with smaller scales due to smoothing, while the sub-problems are more efficiently solved with larger scale parameters.

Third, the number of task pickup-delivery location pairs $|\cJ|$ and the maximum number $K$ of tasks for each driver have a large impact on computational efficiency, because these parameters define the size of the task-chain network (i.e., the number of possible task bundles). Therefore, we test the effect of $\cJ$ with different values of $K$, as shown in panels (c-1) to (c-3). As a result, the CPU time of the N1 mechanism exponentially increases with the increase in $|\cJ|$. When $K = 2$ and $3$, the average CPU time of N1 exceeds 1,000 seconds with large numbers of $|\cJ|$ (despite using the edge-based formulations). Although the CPU time of the FPD mechanism also increases, its increase is approximately linear in $|\cJ|$ and much more modest than that of the N1 mechanism. For instance, our mechanism requires only 16.2 seconds on average when $|\cJ| = 40$ and $K = 3$, which is 330 times more efficient than the N1 mechanism with 5,357.0 seconds on average. Note that the CPU time of the shippers' sub-problems decreases as $|\cJ|$ increases because the number of shippers per shipper group decreases.

%A=B=50000 - macmini
%3.3 seconds for FPD, 897.8 seconds (272 times)
%master 1.2; sub_drv 1.7; sub_shp 0.4

%A=B=50000 - macbookpro
%4.5 seconds for FPD, 1044.6 seconds (232 times)
%master 1.4; sub_drv 2.4; sub_shp 0.7

\begin{figure}[t]
    \centering
    \includegraphics[width=\textwidth]{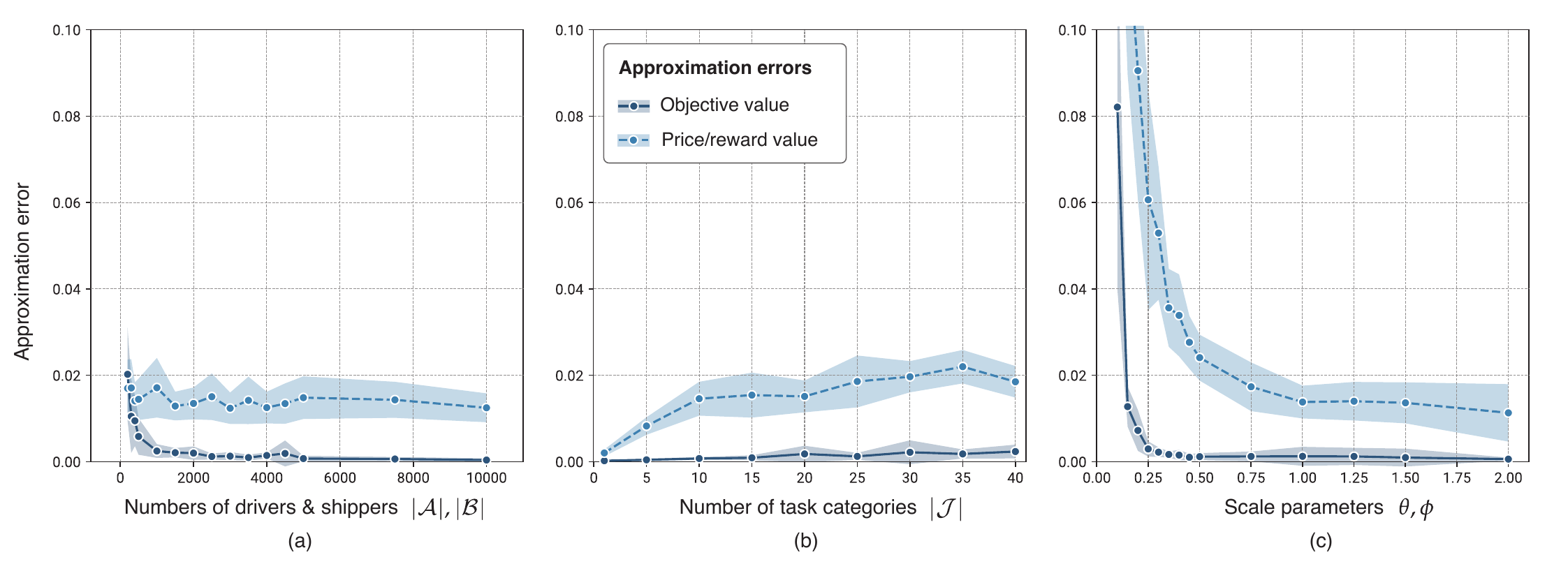}
    \caption{Approximation errors.}
    \label{fig:error}
\end{figure}

\paragraph{Approximation errors.}
\cref{fig:error} shows the approximation errors of the FPD mechanism, with the N1 mechanism considered the ground-truth reference. Overall, the FPD mechanism approximates the N1 mechanism with high accuracy. In most cases the approximation error of the objective is less than 0.5\%, and that of the price/reward value is below 2\%. 

We then discuss the effects of different parameters. Firstly, panel (a) shows the changes in the approximation errors according to the numbers of drivers and shippers (assumed $|\cA| = |\cB|$). The objective error decreases as an increase in $|\cA|$ and $|\cB|$ when these numbers are small. This is because the continuous distributions better approximate the true perceived costs with larger numbers of drivers and shippers. The error becomes 0.3\% on average at $|\cA| = |\cB| = 2000$, and it continues slightly decreasing as $|\cA|$ and $|\cB|$ get larger. Yet, the FPD mechanism achieves the objective error of 2.3\% even when $|\cA| = |\cB| = 200$, the smallest value. The price/reward error does not change with the sizes of $|\cA|$ and $|\cB|$, but the error is always below 2\%. Secondly, panel (b) shows the changes in the approximation errors with the number of task pickup-delivery location pairs $|\cJ|$. Because the number of shippers per OD pair decreases as $|\cJ|$ increases, the continuous approximation gets less accurate; however, the change in the errors is very slight, and the objective error is 0.25\% even when $|\cJ| = 40$.

Finally, panel (c) analyzes the effect of scale parameters $\theta$ and $\phi$ on approximation errors. As larger values of $\theta$ and $\phi$ indicate less heterogeneity in shippers' and drivers' preferences, the FPD mechanism better approximates the socially optimal matching pattern when the scales are larger than certain values. When the scales are small, the approximation errors become large, indicating the need to collect more private cost information through auctions. Still, the objective error is below 2\% when $\theta = \phi = 0.15$ or larger. 

Overall, these results demonstrate the high accuracy of the proposed approach for approximating the socially optimal matching patterns.

\begin{figure}[t]
    \centering
    \includegraphics[width=\textwidth]{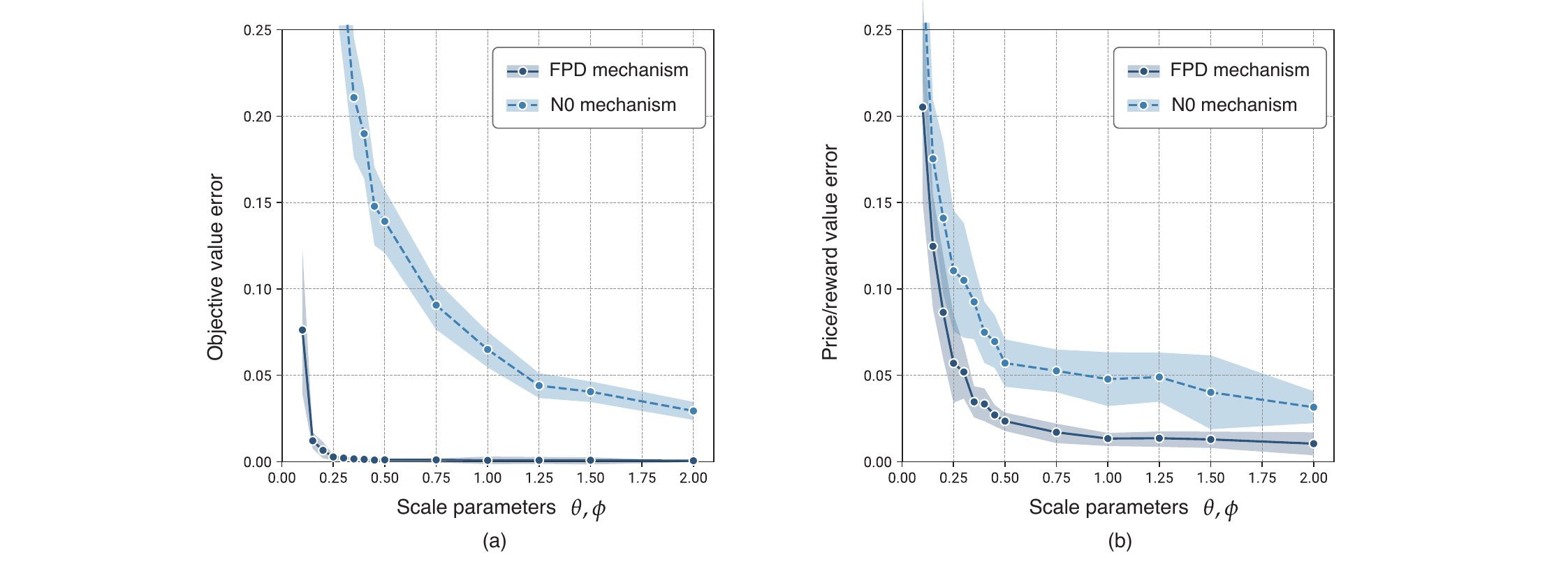}
    \caption{Approximation error comparison between FPD and N0 mechanisms.}
    \label{fig:theta_N0}
\end{figure}

\paragraph{Comparison with the N0 mechanism.}
We finally compare the proposed FPD mechanism with N0 mechanism in terms of approximation errors, to see the impact of capturing heterogeneity. \cref{fig:theta_N0} shows the results: (a) objective errors and (b) price/reward errors, where those of the FPD mechanism are the same as in \cref{fig:error}(c). The approximation errors of N0 mechanism highly depends on the level of heterogeneity in the preferences of drivers and shippers. The N0 mechanism achieves errors less than 5\% when $\theta = \phi = 1.25$ or larger, i.e., the behavior of shippers and drivers is relatively deterministic. However, the erorr become large with a high level of heterogeneity, and both objective and price/reward errors are above 25\% when $\theta = \phi = 0.1$. In contrast, the FPD approach achieves highly accurate solutions even with a high level of heterogeneity, by utilizing the continuous distributions of perceived costs of individual drivers and shippers.

\section{Concluding remarks}\label{sec:conclusion}
This paper presents a general formulation of a two-sided CSD matching problem and proposes a set of methodologies to solve the problem. We first reveal that the FPD approach of \cite{akamatsu_oyama_csd_2024} can be extended to our general formulation. This approach not only addresses computational complexity, but also allows us to collect the information on privately perceived costs for each shipper/driver through VCG auctions. We further transform the master (task partition) problem into a TAP based on a task-chain network representation, which overcomes the difficulty of enumerating task bundles. Then, we formulate the dual problem of the TAP. This enables us to develop an efficient AGD method that only requires a price/reward pattern at market equilibrium, which drastically reduces memory space and computational cost. The numerical experiments demonstrate that our approach solves the general and complex CSD matching problem with high speed ($\sim$ 700 times faster) and accuracy (mostly within 0.5\% errors).

% A limitation of this study is the assumption of task-bundling, where a driver first completes (drops off) a task and then picks up the next task. Although this assumption enables task-chain network representation, our problem is a restricted case of a pickup-delivery routing problem. Therefore, the integration of routing problems is a potential future work of this study. Note that the task-chain network representation can be useful in case the pick-up locations of tasks are consolidated to some depots as assumed in many CSD studies \citep[e.g.,][]{archetti2016vehicle, arslan2019crowdsourced, mancini2022bundle, wang2023joint}.

Finally, we remark that our approach is tested with an MNL model but is compatible with other random/perturbed utility models, including the network multivriate extreme value model \citep{oyama2022markovian}. In addition, the FPD approach and the TAP transformation can also be useful in the ride-sharing context where ordinary drivers serve multiple passengers during their trips through a platform. An extension/application of our approach to a more general context is an important next step of this study.

\appendix
\section{List of notation}\label{app:notation}
Table 1 lists notations frequently used in this paper. %\ref{table:notation}

\begin{table}[t]
\centering
\footnotesize
\label{table:notation}
\caption{Notation in this paper}
\begin{tabular}{lll}
\hline
Category & Symbol & Description  \\ \hline
Sets 
& $\cA$ & Set of drivers \\
& $\cA_{(t,w)}$ & Set of drivers who travel between OD pair $w$ in time-window $t$ \\ 
& $\cB$ & Set of shippers \\
& $\cB_j$ & Set of shippers who want to ship parcel between pickup-delivery location pair $j$ \\ 
& $\cW$ & Set of OD pairs of drivers\\
& $\cJ$ & Set of pickup-delivery location pairs of delivery tasks\\
& $\cT$ ($\tilde{\cT}$) & Set of time windows (with opt-out $0$)\\
& $\cR$ ($\tilde{\cR}$) & Set of task bundles (with opt-out $\emptyset$)\\
& $\cE$ & Set of edges in the task-chain network\\
% & $\cO$ & Set of origins of drivers $\subseteq \cN$\\
% & $\cD$ & Set of destinations of drivers $\subseteq \cN$ \\
% & $\cR$ & Set of pickup nodes of delivery tasks $\subseteq \cN$ \\
% & $\cS$ & Set of destination nodes of delivery tasks $\subseteq \cN$ \\ 
\hline
Parameters 
& $N^D_{(t,w)}$ & Number of drivers in $\cA_{(t,w)}$   \\
& $N^S_j$ & Number of shippers in $\cB_j$ \\
& $K$ & Maximum number of tasks each driver can perform \\
& $\theta$ & Scale of type I EV distribution for shippers' perceived costs \\ 
& $\phi$ & Scale of type I EV distribution for drivers' perceived costs\\ 
\hline
Functions 
&$z$ & Objective function of [SO] \\
&$z^S_j$ & Objective function of shipper sub problem [SO/Sub-Shipper($j$)] \\
&$z^D_{(t,w)}$ & Objective function of driver sub problem [SO/Sub-Driver($t,w$)] \\
&$\Pi^S_j$ & Surplus function of shippers in $\cB_j$ \\
&$\Pi^D_{(t,w)}$, $\pi^D_{(t,w)}$ & Surplus function of drivers in $\cA_{(t,w)}$ \\
&$\cH^S_j$, $\hat{\cH}^S_j$ & Perturbation function of shippers in $\cB_j$ \\
&$\cH^D_{(t,w)}$, $\hat{\cH}^D_{(t,w)}$ & Perturbation function of drivers in $\cA_{(t,w)}$ \\
\hline
Variables 
&$n^b_t$ & 0--1 variable that takes $1$ when shipper $b$ chooses time window $t$ and 0 otherwise \\
&$m^a_r$ & 0--1 variable that takes $1$ when driver $a$ is matched with task bundle $r$ and 0 otherwise \\
&$m^a_{k,ij}$ & 0--1 variable that takes $1$ when driver $a$ is matched with edge $ij$ for $k$-th choice and 0 otherwise \\
& $y_{j,t}$ & Number of shipping permits of delivery time window $t$ issued for shipper group $\cB_j$  \\ 
% & $q^{(t,j)}_w$ & Number of tasks $j$ to be delivered in time window $t$ allocated to driver set $\cA_{(t,w)}$ \\ 
& $f_{(t,w),r}$ & Number of task bundles of $r$ allocated to driver group $\cA_{(t,w)}$\\
& $x^{(t,w)}_{k,ij}$ & Number of drivers in $\cA_{(t,w)}$ to be matched with edge $ij$ for $k$-th choice \\
& $p_{(t,j)}$ & Price/reward for shipping task $j$ in time window $t$ \\
% & $\overline{c}_{rs}$ & Operational cost of a dedicated vehicle to conduct a unit of delivery task $(r,s)$  \\
& $c^b_t$ & Perceived cost for shipper $b$ shipping task in time window $t$ \\
& $C^S_{j,t}$ & Deterministic cost for shippers shipping task $j$ in time window $t$ \\
& $c^a_r$ & Perceived cost for driver $a$ taking detour to perform task bundle $r$ \\ 
& $C^D_{(t,w),r}$ & Deterministic detour cost for drivers who travel between $w$ in $t$ performing task bundle $r$ \\ 
& $\delta_{r,j}$ & Number of tasks between pickup-delivery location pair $j$ included in task bundle $r$\\
& $\delta^r_{k,ij}$ & 0--1 variable that takes 1 if path $r$ uses edge $(i,j)$ at $k$ on the task-chain network and 0 otherwise.\\
% & $\varepsilon^{rs}_{a}$ & Random term of private utility for driver $a$ making delivery task $(r, s)$ \\
% & $t_{kl}$ & Travel time of link $kl\in\cL$, or shortest path travel time between node pair $(k, l)$, $k,l\in\cN$ \\ 
% & $\tau_{ij}$ & Weight of link $(i,j)\in\cL_v$ \\ 
% & $x^{od}_{ij}$ & Flow on link $(i,j)$ whose OD pair is $(o,d)$ \\ 
\hline
% & $\varepsilon$ & Convergence parameter for \Cref{alg:main} \\ \hline
\end{tabular}
\end{table}

\section{Proof of \Cref{prop:so_a_p}}\label{app:proof_so_a}
\begin{proof}
    We apply Lemma 4.2 (and Proposition 4.3) of \cite{akamatsu_oyama_csd_2024} to each of [SO/Sub-Shipper($j$)] and [SO/Sub-Driver($t,w$)] and obtain the objective value functions as follows:
    \begin{align}
        &z^\star_{j}(\mathbf{y}) \approx \sum_{t\in\tilde{\cT}} C^S_{j,t} y_{j,t} - \hat{\cH}^S_j(\mathbf{y}_j) \quad \forall j \in \cJ \label{eq:apprx_opt_val_ship}\\
        &z^\star_{(t,w)}(\mathbf{f}) \approx \sum_{r\in\tilde{\cR}} C^D_{(t,w),r} f_{(t,w),r} - \hat{\cH}^D_{(t,w)}(\mathbf{f}_{(t,w)}) \quad \forall (t,w) \in \cT \times \cW\label{eq:apprx_opt_val_drv}
    \end{align}
    By summing up these over all tasks $j \in \cJ$ and all OD and time-window pairs $(t,w) \in \cT \times \cW$, we obtain the fluidly-approximated master problem [SO/Master-A].
\end{proof}

\section{Proof of \Cref{lemma:so_d}}\label{app:proof_so_d}
\begin{proof}
    The Lagrangian function $\cL$ of [SO/Master-A] is:
    \begin{align}
        \cL(\mathbf{y}, \mathbf{f}, \bm{\lambda}, \bm{\mu}, \bm{\rho}) 
        % &\coloneqq Z(\Theta) 
        % + \sum_{j\in\cJ}\sum_{t\in\cT}\omega_{(t,j)}(y_{j,t} - \sum_{w \in \cW}q^{(t,j)}_{w}) 
        % + \sum_{j \in \cJ} \sum_{t\in\cT} \sum_{w \in \cW} \lambda^{w}_{(t,j)} (q^{(t,j)}_{w} - \sum_{r \in \cR} \delta_{r,j} f_{(t,w),r}) \nonumber\\
        % &~~~~~ + \sum_{j \in \cJ} \mu_{j}(N^S_j - \sum_{t \in \tilde{\cT}} y_{j,t})
        % + \sum_{t\in\cT} \sum_{w \in \cW} \rho^{(t,w)} (N^D_{(t,w)} - \sum_{r \in \cR} f_{(t,w),r}) \nonumber\\
        &= \sum_{j \in \cJ} \qty( \sum_{t \in \cT} y_{j,t} (C^S_{j,t} - \mu_j + \lambda_{(t,j)}) - \hat{\cH}^S_j(\mathbf{y}_j) ) \nonumber\\
        &~~~+ \sum_{t\in\cT} \sum_{w \in \cW} \qty(
        \sum_{r \in \tilde{\cR}} f_{(t,w),r} (C^D_{(t,w),r} - \sum_{j \in \cJ} \delta_{r,j}\lambda_{(t,j)} - \rho_{(t,w)}) - \hat{\cH}^D_{(t,w)}(\mathbf{f}_{(t,w)}) ) \nonumber\\
        &~~~+ \sum_{j \in \cJ} \mu_{j}N^S_j + \sum_{t\in\cT} \sum_{w \in \cW} \rho_{(t,w)} N^D_{(t,w)}
    \end{align}
    where $\bm{\lambda}, \bm{\mu}, \bm{\rho}$ are the Lagrangian multipliers with respect to \eqref{eq:pathflow_leq_assign}, \eqref{eq:N_ship_consv}, \eqref{eq:N_drv_consv}, respectively.
    % Rearranging the terms with respect to $\mathbf{y}$ and $\mathbf{f}$ yields:

    From the optimality conditions with respect to $\mathbf{y}, \mathbf{f},\bm{\lambda}$, we obtain:
    \begin{align}
        \begin{dcases}
            y^\star_{j,t} (C^S_{j,t} - \mu_{j} + \omega_{(t,j)} - \hat{\cH}^\prime_{j,t}) = 0, \, y^\star_{j,t} \geq 0, \, C^S_{j,t} - \mu_{j} + \lambda_{(t,j)} - \hat{\cH}^\prime_{j,t} \geq 0\\
            f^{(t,w)\star}_r (C^D_{(t,w),r} - \sum_{j \in \cJ} \delta_{r,j}\lambda_{(t,j)} - \rho_{(t,w)} - \hat{\cH}^\prime_{r}) = 0, \, f^{(t,w)\star}_r \geq 0, \, C^D_{(t,w),r} - \sum_{j \in \cJ} \delta_{r,j}\lambda^\star_{(t,j)} - \rho_{(t,w)} - \hat{\cH}^\prime_{r} \geq 0\\
            \lambda^\star_{(t,j)} (\sum_{w\in\cW}\sum_{r\in\tilde{\cR}}f^\star_{(t,w),r} - y^\star_{j,t}) = 0, \, \lambda^\star_{(t,j)} \geq 0 ,\,  \sum_{w\in\cW}\sum_{r\in\tilde{\cR}}f^\star_{(t,w),r} - y^\star_{j,t} \geq 0.
        \end{dcases}
    \end{align}
    The first and second conditions are the random/perturbed utility maximization principles of shippers and drivers, and the third is the market clearing condition. Thus, we can interpret $\mu_j$ and $\rho_{(t,w)}$ as the expected minimum costs $\Pi^S_j$ and $\Pi^D_{(t,w)}$, and $\lambda_{(t,j)} = p_{(t,j)}$. Moreover, as the perturbed functions $\hat{\cH}^S_j$ and $\hat{\cH}^D_{(t,w)}$ are the convex conjugate of $\Pi^S_j$ and $\Pi^D_{(t,w)}$,
    \begin{align}
        &N^S_j  \Pi^S_j(\mathbf{p}) = \sum_{t \in \cT} y^\star_{j,t} (C^S_{j,t} + p_{(t,j)}) - \hat{\cH}^S_j(y^\star)\\
        %\sum_{t \in \cT} y^\star_{j,t} (C^S_{j,t} - \mu_j + \omega_{(t,j)}) - \hat{\cH}^S_j(y^\star) = \qty(\sum_{t \in \cT} y^\star_{j,t} (C^S_{j,t} + \omega_{(t,j)}) - \hat{\cH}^S_j(y^\star)) - N^S_j  \mu_j = 0\\
        &N^D_{(t,w)} \Pi^D_{(t,w)}(\mathbf{p}) = \sum_{r \in \tilde{\cR}} f^\star_r (C^D_{(t,w),r} - \sum_{j \in \cJ} \delta_{r,j}p_{(t,j)}) - \hat{\cH}^D_{(t,w)}(\mathbf{f}^\star)
        % \sum_{r \in \cR} f^\star_r (c_r - \sum_{j \in \cJ} \delta_{r,j}\lambda_{(t,j)} - \rho_{(t,w)}) - \hat{\cH}^D_{(t,w)}(\mathbf{f}^\star) = \qty(\sum_{r \in \cR} f^\star_r (c_r - \sum_{j \in \cJ} \delta_{r,j}\lambda_{(t,j)}) - \hat{\cH}^D_{(t,w)}(\mathbf{f}^\star)) - N^D_{(t,w)}\rho_{(t,w)} = 0
    \end{align}
    As a result, the Lagrangian dual problem is
    \begin{align}
    \max_{\bm{\lambda}, \bm{\mu}, \bm{\rho}} \min_{\mathbf{y}, \mathbf{f}} \cL(\mathbf{y}, \mathbf{f}, \bm{\lambda}, \bm{\mu}, \bm{\rho}) 
    = \max_{\mathbf{p}} \cL(\mathbf{y}^\star, \mathbf{f}^\star, \mathbf{p}) = 
    \sum_{j \in \cJ} N^S_j \Pi^S_j(\mathbf{p})  +
    \sum_{t\in\cT}\sum_{w \in \cW} N^D_{(t,w)} \Pi^D_{(t,w)}(\mathbf{p}) 
    \label{eq:L_dual}
    \end{align}
    Moreover, $\Pi^D_{(t,w)}$ is equivalent to $\pi^{(t,w)}_{(0,o)}$, the expected minimum cost from the source state $(0,o)$ to the target state $(K+1, d)$ in the task-chain network of \Cref{fig:driver_net}.
\end{proof}

\section{Proof of \Cref{prop:VCG_ship}}\label{sec:proof_VCG_ship}
\begin{proof}
    Under the payment \eqref{eq:vcg_reward_ship}, the utility that shipper $b$ gains is
    \begin{align}
        \label{eq:vcg_shipper_util}
        u_b(\mathbf{s}_b, \mathbf{s}_{-b}) &= -\sum_{t\in\tilde{\cT}} n^{b\star}_t(\mathbf{s}_b, \mathbf{s}_{-b}) c^b_t - p_b(\mathbf{s}_b, \mathbf{s}_{-b}) \nonumber\\
        &= \qty(-\sum_{t\in\tilde{\cT}} n^{b\star}_t(\mathbf{s}_b, \mathbf{s}_{-b}) c^b_t + \sum_{b'\neq b}\sum_{t\in\tilde{\cT}} n^{b'\star}_t(\mathbf{s}_b, \mathbf{s}_{-b}) s^{b'}_t) - \max_{\mathbf{n}_{-b}} z^S_{j,-b}(\mathbf{n}_{-b}|\mathbf{s}_{-b},\mathbf{y}) %\in \cB_j\setminus\{b\}
    \end{align}
    The final term is independent of $b$, and shipper $b$ would choose $\mathbf{s}_b$ to maximize 
    \begin{align}\nonumber
        -\sum_{t\in\tilde{\cT}} n^{b\star}_t(\mathbf{s}_b, \mathbf{s}_{-b}) c^b_t + \sum_{b'\neq b}\sum_{t\in\tilde{\cT}} n^{b'\star}_t(\mathbf{s}_b, \mathbf{s}_{-b}) s^{b'}_t.
    \end{align}
    Considering the assignment rule, $n^{b\star}_t(-\mathbf{c}_b, \mathbf{s}_{-b})$ is the matching pattern that maximizes this. Therefore, for any $\mathbf{s}_b \neq -\mathbf{c}_b$, the following holds:
    \begin{align}
        u_b(-\mathbf{c}_b, \mathbf{s}_{-b}) \geq u_b(\mathbf{s}_b, \mathbf{s}_{-b}),
    \end{align}
    that is, the optimal behavior for $b$ to maximize their utility is to announce the true valuation $\mathbf{s}_b = -\mathbf{c}_b$, which results in maximizing the social surplus. Thus, [VCG/Shipper] is \textit{truthful} and \textit{socially efficient} \citep[see also e.g.,][]{cramton2006}. 
\end{proof}

\section{Proof of \Cref{prop:VCG_drv}}
\begin{proof}\label{sec:proof_VCG_drv}
    Under the reward scheme \eqref{eq:vcg_reward_drv_edge}, the utility of driver $a$ is
    \begin{align}
        \label{eq:vcg_driver_util}
        u_a(\tilde{\mathbf{c}}_a, \tilde{\mathbf{c}}_{-a}) &= w_{a}(\tilde{\mathbf{c}}_a, \tilde{\mathbf{c}}_{-a}) - \sum_{k=0}^{K} \sum_{ij\in\cE}
            m^{a\star}_{k,ij}(\tilde{\mathbf{c}}_a, \tilde{\mathbf{c}}_{-a}) c^{a}_{ij} \nonumber\\
        &= \min_{\mathbf{m}_{-a}} z^D_{-a}(\mathbf{m}_{-a}|\tilde{\mathbf{c}}_{-a},\mathbf{x}) 
            - \qty(\sum_{k=0}^{K} \sum_{ij\in\cE}
            m^{a\star}_{k,ij}(\tilde{\mathbf{c}}_a, \tilde{\mathbf{c}}_{-a}) c^{a}_{ij} +\sum_{a'\neq a}\sum_{k=0}^{K} \sum_{ij\in\cE}
            m^{a'\star}_{k,ij}(\tilde{\mathbf{c}}_a, \tilde{\mathbf{c}}_{-a}) \tilde{c}^{a'}_{ij}). %\in\cA_{(t,w)}\setminus\{a\}
    \end{align}
    Because the first term of the RHS is independent of $a$, driver $a$ would choose bid $\tilde{\mathbf{c}}_a$ to minimize
    \begin{align}\nonumber
        \sum_{k=0}^{K} \sum_{ij\in\cE}
            m^{a\star}_{k,ij}(\tilde{\mathbf{c}}_a, \tilde{\mathbf{c}}_{-a}) c^{a}_{ij} +\sum_{a'\neq a}\sum_{k=0}^{K} \sum_{ij\in\cE}
            m^{a'\star}_{k,ij}(\tilde{\mathbf{c}}_a, \tilde{\mathbf{c}}_{-a}) \tilde{c}^{a'}_{ij}
    \end{align}
    Considering the assignment rule of [VCG/Driver] (in Step 2), $m^{a'\star}_{k,ij}(\mathbf{c}_a, \tilde{\mathbf{c}}_{-a})$ would minimize this. Therefore, for any $\tilde{\mathbf{c}}_a \neq \mathbf{c}_a$, the following holds:
    \begin{align}
        u_a(\mathbf{c}_a, \tilde{\mathbf{c}}_{-a}) \geq u_a(\tilde{\mathbf{c}}_a, \tilde{\mathbf{c}}_{-a}),
    \end{align}
    that is, announcing true perceived cost $\mathbf{c}_a$ by bidding is a (weakly) dominant strategy for each driver $a$, and in this case, a socially optimal matching is achieved. This proves that [VCG/Driver] is \textit{truthful} and \textit{efficient}.
    
    % Using this function, we first prove the truthfulness by contradiction.
    % Assume that the following inequality holds:
    % \begin{align}
    %     \label{eq:contradiction}
    %     u_a(\mathbf{c}_a, \tilde{\mathbf{c}}_{-a}) < u_a(\tilde{\mathbf{c}}_a, \tilde{\mathbf{c}}_{-a}).
    % \end{align}
    % By subsutituting \eqref{eq:vcg_driver_util} into \eqref{eq:contradiction},
    % \begin{align}
    %     \sum_{k=0}^{K} \sum_{ij\in\cE}
    %         m^{a\star}_{k,ij}(\mathbf{c}_b, \tilde{\mathbf{c}}_{-b}) c^{a}_{ij} +\sum_{a'\neq a}\sum_{k=0}^{K} \sum_{ij\in\cE}
    %         m^{a'\star}_{k,ij}(\mathbf{c}_b, \tilde{\mathbf{c}}_{-b}) \tilde{c}^{a'}_{ij} >
    %     \sum_{k=0}^{K} \sum_{ij\in\cE}
    %         m^{a\star}_{k,ij}(\tilde{\mathbf{c}}_a, \tilde{\mathbf{c}}_{-a}) c^{a}_{ij} +\sum_{a'\neq a}\sum_{k=0}^{K} \sum_{ij\in\cE}
    %         m^{a'\star}_{k,ij}(\tilde{\mathbf{c}}_a, \tilde{\mathbf{c}}_{-a}) \tilde{c}^{a'}_{ij}
    % \end{align}
    % and thus
    % \begin{align}
    %     z^{D\star}(\mathbf{c}_b, \tilde{\mathbf{c}}_{-b})
    %      >
    %     z^{D\star}(\tilde{\mathbf{c}}_b, \tilde{\mathbf{c}}_{-b})
    %     +
    %     \sum_{k=0}^{K} \sum_{ij\in\cE}
    %         m^{a\star}_{k,ij}(\tilde{\mathbf{c}}_a, \tilde{\mathbf{c}}_{-a}) (c^{a}_{ij} - \tilde{c}^{a}_{ij})
    % \end{align}
\end{proof}

\section{Edge-based formulation of [SO]}\label{app:so_edge}
Due to the difficulty of enumerating task bundles, in the numerical experiments of \Cref{sec:experiment} we solve [SO] using an edge-based formulation on the task-chain network, as follows.
\begin{align}
    [\text{SO-Edge}]\qquad
    \min_{\mathbf{m},\mathbf{n}}\quad&
    z(\mathbf{m},\mathbf{n}) \equiv \sum_{j\in\cJ}\qty(\sum_{b\in\cB_{j}} \sum_{t\in\tilde{\cT}} n^b_t c^b_t) +
    \sum_{t\in\cT}\sum_{w \in \cW}\qty(\sum_{a\in\cA_{(t,w)}}
    \sum_{k=0}^{K} \sum_{ij\in\cE}
    m^a_{k,ij} c^a_{ij})
    \label{eq:SO_obj}\\
    \subto\quad& 
    \sum_{t\in\tilde{\cT}} n^b_t = 1 \qquad \forall b \in \cB, \tag{\ref{eq:indv_shipper_consv_all}}\\
    &n^b_t\in\{0,1\} \qquad \forall t \in\tilde{\cT}, \, \forall b \in \cB, \tag{\ref{eq:n_binary_all}}\\
    &\sum_{i\in\tilde{\cJ}} m^a_{k-1,ij} - \sum_{l\in\tilde{\cJ}} m^a_{k,jl} = \eta_{k,j}  \quad \forall k \in \{0,\ldots,K+1\}, \, \forall j \in \cJ,\label{eq:m_edge_flow_consv}\\
    & m^a_{k,ij} \in \{0,1\}  \quad \forall k \in \{0,\ldots,K+1\}, \, \forall ij \in \cE, \, \forall a \in \cA,\label{eq:m_edge_binary_all}\\
    & \sum_{b\in\cB_{j}}n^b_t \leq \sum_{w \in \cW}\sum_{a\in \cA_{(t,w)}}\sum_{k=0}^{K} \sum_{i\in\cJ} m^a_{k,ij} \qquad \forall j \in \cJ, \,  \forall t \in \cT
\end{align}

\bibliographystyle{elsarticle-harv} %plainnat
\bibliography{cite} 

\end{document}